\documentclass[a4paper,12pt]{article}
\pdfoutput=1
\usepackage{authblk}
\usepackage{amsmath}
\numberwithin{equation}{section}
\usepackage{amssymb}
\usepackage{amsthm}
\usepackage{bm}
\usepackage{color}
\usepackage{array}
\usepackage{geometry}
\usepackage{hyperref}
\hypersetup{hidelinks}
\usepackage{longtable}
\usepackage{multirow}
\usepackage{makecell}
\usepackage{enumerate}
\usepackage{graphicx}
\usepackage{subfigure}
\usepackage{epstopdf}
\usepackage{embrac}
\usepackage[numbers,sort&compress]{natbib} 
\usepackage{algorithm,algorithmic}
\numberwithin{algorithm}{section}
\newtheorem{thm}{Theorem}
\newtheorem{ass}{Assumption}
\newtheorem{lem}{Lemma}
\newtheorem{rem}{Remark}
\newtheorem{defi}{Definition}
\newtheorem{propo}{Proposition}
\newtheorem{examp}{Example}
\numberwithin{thm}{section}
\numberwithin{lem}{section}
\numberwithin{rem}{section}
\numberwithin{defi}{section}
\numberwithin{propo}{section}
\numberwithin{examp}{section}

\newcommand{\G}{\mathcal{G}}
\newcommand{\V}{\mathcal{V}}
\newcommand{\E}{\mathcal{E}}
\newcommand{\mR}{\mathbb{R}}

\newcommand{\1}{{\bf 1}}
\newcommand{\0}{{\bf 0}}
\newcommand{\proxTi}{{\rm prox}_{\alpha T_i}}
\newcommand{\proxT}{{\rm prox}_{\alpha T}}
\newcommand{\proxtP}{{\rm prox}_{P^{-1}\Phi}\circ P^{-1}Q}
\title{Distributed Proximal-Correction Algorithm for the Sum of Maximal Monotone Operators in Multi-Agent Network}
\date{}
\author[1]{Kai Gong }
\author[1]{Liwei Zhang}
\affil[1]{School of Mathematical Sciences, Dalian University of Technology, Dalian,
	116024, Liaoning, China. \authorcr Email: gk1995\_\_\_@mail.dlut.edu.cn
	\authorcr Email: lwzhang@dlut.edu.cn}
\begin{document}
	\maketitle
	\begin{abstract}
        This paper focuses on a class of inclusion problems of \emph{maximal monotone} operators in a multi-agent network, where each agent is characterized by an operator that is not available to any other agents, but the agents can cooperate by exchanging information with their neighbors according to a given communication topology. All agents aim at finding a common decision vector that is the solution to the sum of agents' operators. This class of problems are motivated by distributed convex optimization with coupled constraints. In this paper, we propose a distributed proximal point method with a cumulative correction term (named Proximal-Correction Algorithm) for this class of inclusion problems of operators. It's proved that Proximal-Correction Algorithm converges for any value of a constant penalty parameter. In order to make Proximal-Correction ALgorithm being computationally implementable for a wide variety of distributed optimization problems, we adopt two inexact criteria for calculating the proximal steps of the algorithm. Under each of these two criteria, the convergence of Proximal-Correction Algorithm can be guaranteed, and the linear convergence rate is established when the stronger one is satisfied. In numerical simulations, both exact and inexact versions of Proximal-Correction Algorithm are executed for a distributed convex optimization problem with coupled constraints. Compared with several alternative algorithms in literature, the exact and inexact versions of Proximal-Correction both exhibit good numerical performance.\\
	   	\vspace{0.5em}
	
	\end{abstract}
	\section{Introduction}
      In this paper, we consider a connected and undirected graph $\G=(\V,\E)$ with vertex set $\V=\{i=1,2,\dots,N\}$, and edge set $\E=\{(i,j):i,j\in\V\}$, and intend to solve the following inclusion problem of the sum of $N$ maximal monotone operators over $\G$:
      \begin{equation}\label{P}
      	\text{determining} \quad z\in\mR^n \quad \text{such that} \quad \0\in\tilde{T}(z)\triangleq\sum_{i=1}^NT_i(z),
      \end{equation}
      where each $T_i:\mR^n\rightrightarrows\mR^n, i=1,2,\dots,N$, is a maximal monotone operator.
	  In problem \eqref{P}, for each $i=1,2,\dots,N$, agent $i$ has a decision variable $z_i$, the local copy of the global decision variable $z$, and has access only to it's local operator $T_i$, but it can exchange its local decision vectors with its neighbors according to $\G$. All agents cooperate to reach an aligned decision vector that is the solution to problem \eqref{P}.
	  
	  The relation $\0\in\sum_{i=1}^NT_i(z)$ can be viewed as an extension of the first-order optimality conditions, hence it covers a wider range of distributed convex optimization problems. For example, if $T_i$ is the gradient $\nabla f_i$ (or the subdifferential $\partial f_i$) of a convex local objective function $f_i$, obviously, in this case, \eqref{P} is to solve a distributed unconstrained convex optimization problem. Also, other types of distributed convex optimization problems, including the type where the coupled inequality and coupled affine equality constrains are present, can be formulated as special cases of problem \eqref{P}. The details would be given in Section 2. Next, we would like to review some related works on the important applications of problem \eqref{P}, i.e. distributed convex optimization problems. 
	  
	  For unconstrained distributed convex optimization problem, there are many classical algorithms such as DGD \cite{Nedic09}, EXTRA \cite{EXTRA}, DIGing \cite{DIGing}, NEXT \cite{NEXT}, etc. On the other hand, different types of distributed constrained optimization problems have been deeply studied and efficient algorithms for such problems have been developed. For the distributed convex optimization with local constraint sets, when the local objective functions have Lipschitz continuous gradients, both PSCOA \cite{YuWW21} and PG-EXTRA \cite{PG-EXTRA} converge with a rate of $O(1/k)$, NIDS \cite{NIDS} achieves a faster rate of $o(1/k)$. Furthermore, if the local objective functions are also strongly convex, the gradient-tracking-based distributed method \cite{LiuHZ22} reaches a linear rate. In addition, distributed optimization problems with coupled inequality and coupled affine equality constraints are encountered in a wide range of applications such as power systems, plug-in electric vehicles charging problems, and optimal resource allocation problems, and have been investigated intensively in recent years. To name a few, each of the distributed proximal primal-dual (DPPD) method \cite{LiXX21}, the gradient-based distributed primal-dual method \cite{LiangS21}, the proximal dual consensus ADMM (PDC-ADMM) method \cite{Chang16}, and the the distributed integrated primal-dual proximal (IPLUX) method \cite{WuXY22} is guaranteed to converge to an optimal solution. Among these algorithms, PDC-ADMM and IPLUX have a convergence rate of $O(1/k)$.
	  
	  Maximal monotone operators play important roles in convex analysis and have been studied extensively.
	  To find a solution to a maximal monotone operator in centralized settings, the proximal point algorithm (PPA) \cite{Rockafellar76} is an efficient method and has been employed to construct various distributed algorithms \cite{LiXX19,LiXX21,Falsone22,WuXY22,XuJM21}.  In addition, the work in \cite{Rockafellar76} adopted two inexact rules for calculating the proximal points, which enables PPA to be computationally implementable for a wide variety of problems. Such a relaxation  therefore can be taken by distributed algorithms based on PPA. With such observations, we propose the distributed Proximal-Correction Algorithm and its inexact form, which exhibit the following features:
	  \begin{enumerate}[\hspace{1em}(1)]
	  	\item { Proximal-Correction Algorithm can be employed for a wide variety of distributed convex optimization problems, such as the listed examples in section 2;}
  	    \item {local objective and constraint functions are not required to be differentiable, Lipschitz continuous, or smooth, but only convex and lower semicontinuous;}
  	    \item {convergence of Proximal-Correction Algorithm is guaranteed for any value of a constant penalty parameter. In addition, the proximal steps of Proximal-Correction Algorithm can be approximately calculated under two criteria, one of which guarantees the convergence to an optimal solution and the stronger another provides the linear convergence rate.}
	  \end{enumerate}
      Algorithms for distributed optimization problems with coupled constraints generally require to solve minimization subproblems. To the best of our knowledge, none of the previously mentioned distributed algorithms in literature adopted inexact computations of subproblems. In this article we take inexact rules for solving subproblems in the distributed setting.

	 The rest of the article is divided into five sections. Section 2 contains some necessary preliminaries and exapmles of distributed optimization problems included in problem \eqref{P}. Section 3 first constructs a distributed algorithm by employing the proximal point method, and then considers two criteria to approximate the proximal steps in this algorithm. Section 4 develops the convergence analysis of the proposed algorithms including the exact and the inexact proximal steps, and presents the main theoretical results. Section 5 demonstrates the effectiveness of these two algorithms through numerical simulations and comparison with several other alternative algorithms. Section 6 makes conclusions.
	 
	 \section{Preliminaries}
	 In this section, we first introduce some preliminaries, including required concepts and notations. Then we give several examples of the applications of problem \eqref{P} in distributed convex optimization and present the required standard assumptions.
	 \subsection{Notions and Notations}
	 The $n$-dimensional Euclidean space is denoted by $\mR^n$. We write $\1$ and $\0$ as a vector with all ones and all zeros elements, respectively. For a nonempty subset $\Omega$ in $\mR^n$, the  interior of $\Omega$ is written as ${\rm int}(\Omega)$. $I_\Omega(x)$ denotes the indicator function of $\Omega$, i.e., whenever $x\in\Omega$, the function
	 has value $0$, otherwise it is $+\infty$. We write $K^\circ$ as the polar of a cone $K$. For a matrix $W$, denote by $W^\top$ the transpose of $W$. $W\succ\0 \ (W\succcurlyeq\0)$ means that $W$ is symmetric positive definite (semidefinite), and $B\succcurlyeq A$ indicates $B-A\succcurlyeq \0$. ${\rm null}\{W\}$ represents the null space, and $\rho(W)$ is the spectral radius of $W$. The Kronecker product is denoted by $\otimes$. ${\rm span}\{x\}$ is the linear space generated by vector $x$. $\partial f(x)$ is the subdifferential of a function $f$ at $x$. Let $T: H\rightrightarrows H$ be a set-valued mapping, where $H$ is a Hilbert space. We write $T^{-1}$ as the inverse mapping of $T$, the set $D(T)=\{x| T(x)\neq\emptyset\}$ as its effective domain, and $G(T)=\left\{(x,u)\in(H\times H)| u\in T(x)\right\}$ as its graph. Suppose $A$ is a symmetric positive definite matrix, we define the $A$-matrix norm of $x$ as
	 $\Vert x\Vert_A\triangleq\sqrt{\langle x, Ax\rangle}$, and $A$ is ignored if A is the identity matrix $I$. The projection operator is abbreviated as ${\rm Proj}$.
	
	 \begin{defi}[maximal monotone operator]
	 	Let $H$ be a real Hilbert space with inner product $\langle\cdot,\cdot\rangle$, A multifunction $T: H\rightrightarrows H$ is said to be a monotone operator if
	 	\begin{equation*}
	 		\langle x-y, u-v\rangle\ge0,\quad\text{whenever}\quad u\in T(x), v\in T(y).
	 	\end{equation*}
	 	It is said to be maximal monotone if, in addition, the graph $G(T)$
	 	is not properly contained in the graph of any other monotone operator $T^\prime: H\rightrightarrows H$.
	 \end{defi}
	 Note that the maximality of a monotone operator $T$ implies that for any pair $(\hat{x},\hat{u})\notin G(T)$, there exists $(x,u)\in G(T)$ with $\langle x-\hat{x}, u-\hat{u}\rangle<0$. This property will be often used in convergence analysis.
	 \begin{defi}(\cite{Rockafellar76})
	 	We shall say that $\Phi^{-1}$, the inverse of a multifunction $\Phi$, is Lipschitz continuous at $\0$ with modulus $a\ge0$, if there is a unique solution $\xi^*$ to $\0\in\Phi(\xi)$ \emph{(i.e. $\Phi^{-1}(0)=\{\xi^*$\})}, and for some $\tau>0$ we have
	 	\begin{equation}\label{ine_Lip}
	 		\Vert \xi-\xi^*\Vert\le a\Vert w\Vert \quad \text{whenever} \ \xi\in\Phi^{-1}(w) \ \text{and} \ 
	 		\Vert w\Vert\le\tau.
	 	\end{equation} 
	 \end{defi}
	 \subsection{Examples and Assumptions}
	 Problem \eqref{P} contains a wide category of distributed convex optimization. We would like to provide  some particular examples of problem \eqref{P} to illustrate the value of this fundamental mathematical model.
	 \begin{examp}
	 	When the operator $T_i$ is the subdifferential $\partial f_i$ of a convex function $f_i:\mR^n\to (-\infty,+\infty], f_i\not\equiv+\infty$, then $T_i$ is maximal monotone. The relation $0\in\sum_{i=1}^NT_i(z)$ means that 
	 	\begin{equation*}
	 	z\in\mathop{\arg\min}\limits_{x\in\mR^n}\sum_{i=1}^Nf_i(x).
	 	\end{equation*}
	 \end{examp}
	 \begin{examp}\label{Examp_2}
	 Consider a distributed convex optimization problem with coupled constraints, i.e.,
	 \begin{equation}\label{DOPCC}
	 	\begin{aligned}
	 		&\min_{x\in\mR^s}f(x)\triangleq\sum_{i=1}^Nf_i(x)\\
	 		&{\rm s.t.}\quad\sum_{i=1}^Ng_i(x)\le\0,
	 		\quad\sum_{i=1}^N A_ix=b,\quad x\in\bigcap_{i=1}^N\Omega_i,
	 	\end{aligned}
	 \end{equation}
    where $f_i:\mR^s\to\mR$ and $g_i:\mR^s\to\mR^p$ are lower semicontinuous and convex functions, $\Omega_i$ is a closed and convex subset in $\mR^s$, for each $i=1,2,\dots,N$. We would like to transform problem \eqref{DOPCC} to a special case of problem \eqref{P} by using dual decomposition.
     Denote 
    \begin{align*}
    	K\triangleq\mR^p_-\times\{\0_q\},\quad K^\circ=\mR^p_+\times\mR^q,\quad
    	h_i(x)\triangleq\left[
    	\begin{array}{c}
    		g_i(x)\\
    		A_ix-b/N
    	\end{array}
    	\right].
    \end{align*}
    The local Lagrange function $L_i:\mR^s\times\mR^m\to (-\infty,+\infty]$ satisfies
    \begin{equation}\label{eq_Lag_Fun}
    	L_i(x,\lambda)=f_i(x)+\lambda^\top h_i(x)+I_{\Omega_i}(x)-I_{K^\circ}(\lambda),\quad i=1, 2,\dots, N.
    \end{equation}
    Then the dual of problem \eqref{DOPCC} can be formulated as 
    \begin{equation}\label{DP}
    	\max_{\lambda\in\mR^m}\min_{x\in\mR^s}\sum_{i=1}^NL_i(x,\lambda).
    \end{equation}
    Under Slater's condition, there is no duality gap between problem \eqref{DOPCC} and \eqref{DP}. We turn to solve the dual problem \eqref{DP} to obtain a solution to \eqref{DOPCC}, i.e.,
    \begin{equation*}
    	\text{determining}\quad (x,\lambda)\in\mR^s\times\mR^m\quad \text{ such that}\quad \0\in\sum_{i=1}^N\partial L_i(x,\lambda).
    \end{equation*}
	 \end{examp}
     \begin{examp}
      If the agents' decision variables in problem \eqref{DOPCC} are not consensual, this problem is considered in \cite{WuXY22,ShuL20,LiangS21}, i.e.,
     	\begin{align*}
     		&\min_{x_1,\dots,x_N}f(x):=\sum_{i=1}^Nf_i(x_i)\\
     		&{\rm s.t.}\quad \sum_{i=1}^Ng_i(x_i)\le\0,\
     		\ \sum_{i=1}^N A_ix_i=b,\quad x_i\in\Omega_i, i=1,2,\dots,N.
     	\end{align*}
     	By the dual decomposition, we obtain its dual
     	\begin{equation*}
     		\max_{\lambda\in\mR^m}\sum_{i=1}^Nd_i(\lambda)
     	\end{equation*} 
     	as a special case of problem \eqref{P} with $T_i=-\partial d_i$, where $d_i(\lambda)=\mathop{\min}\limits_{x_i}L_i(x_i,\lambda)$, and $L_i$ is defined by \eqref{eq_Lag_Fun}.
     \end{examp}
     \begin{rem}
     	In these three Examples, note that we do not require the local objective functions $f_i$'s and the constraint functions $g_i$'s to be differentiable or smooth, they are just lower semicontinuous and convex. Besides, we also do not require the local constraint sets $\Omega_i, i=1, 2,\dots, N$, to be neither bounded nor equal to each other. 
     \end{rem}
     \begin{examp}
     	In a centralized setting, the relation $0\in \tilde{T}(z)$ can be linked to the well-known complementarity problems of mathematical programming \cite{Rockafellar76}, if the operator $\tilde{T}$ satisfies
     	\begin{equation}\label{eq_T_eg}
     		\begin{aligned}
     			\tilde{T}(z)=\left\{
     			\begin{array}{ll}
     				T_0(z)+N_D(z)& \text{if} \quad z\in D,\\
     				\emptyset& \text{if} \quad z\notin D,
     			\end{array}
     			\right.
     		\end{aligned}
     	\end{equation}
     	where the operator $T_0$ is single-valued, and $N_D(z)$ represents the normal cone of a nonempty closed convex subset $D$ at $z$. Naturally, such problems also can be considered in the decentralized environment as each local operator $T_i$ has an expression like \eqref{eq_T_eg}, such an extension may be potentially worthwhile. 
     \end{examp}   
     
	 With the observations of these examples listed above, in the following sections we will design an algorithm to solve the fundamental problem \eqref{P} rather than considering only its special examples.
	  We make the following assumption to guarantee that a solution to \eqref{P} exists.
	 \begin{ass}\label{ass_1}
	 	\hspace{1em}
	 	\begin{enumerate}[(a)]
	 		\item {Each local operator $T_i$ of agent $i$ is maximal monotone, $i=1, 2,\dots, N$.}
	 		\item {$\left(\bigcap_{i=1}^{N-1}D(T_i)\right)\bigcap{\rm int}\left(D(T_N)\right)\neq\emptyset$.}
	 		\item {Problem \eqref{P} admits one solution $z^*$.}
	 	\end{enumerate}
	 \end{ass}
	 Under Assumption \ref{ass_1}, it is known that the sum $\tilde{T}$ is also maximal monotone \cite{Rockafellar70}, then there exists at least one fixed point of its proximal mapping $(I+\alpha \tilde{T})^{-1}$ for any $\alpha>0$ (see \cite{Rockafellar76}), which is exactly a solution to \eqref{P}. In this paper, we consider a connected undirected network and share the same assumptions with \cite{EXTRA} on the mixing matrices of the multi-agent network as follows.
	 \begin{ass}[mixing matrix]\label{ass_2}
	 	Consider a connected network $\G=\left\{\V,\E\right\}$ consisting of a set of agents $\V=\{1, 2, \dots, N\}$ and a set of undirected edges $\E$. The mixing matrices $W=[w_{ij}]\in\mR^{N\times N}$ and $\tilde{W}=[\tilde{w}_{ij}]\in\mR^{N\times N}$ satisfy 
	 	\begin{enumerate}[(a)]
	 		\item {If $i\neq j$ and $(i,j)\notin\E$, then $w_{ij}=\tilde{w}_{ij}=0$.}\label{ass_2a}
	 		\item {$W\1=\1, \1^\top W=\1^\top$,\quad$\tilde{W}\1=\1, \1^\top\tilde{W}=\1^\top$.}\label{ass_2b}
	 		\item {${\rm null}\{\tilde{W}-W\}={\rm span}\{\1\},\quad {\rm null}\{I-\tilde{W}\}\supseteq {\rm span}\{\1\}$.}\label{ass_2c}
	 		\item {$\tilde{W}\succ\0$,\quad $(I+W)/2\succcurlyeq\tilde{W}\succcurlyeq W$.}\label{ass_2d}
	 	\end{enumerate} 
	 \end{ass}
	 It is known that parts \eqref{ass_2b}-\eqref{ass_2d} of Assumption \ref{ass_2} imply 
	 ${\rm null}\{I-W\}={\rm span}\{\1\}$ and the eigenvalues of $W$ lie in $(-1, 1]$, hence the matrix $(I+W)/2$ is positive definite, and we can simply take $\tilde{W}=(I+W)/2$ that is found to be very efficient. There are a few common choices of mixing matrix $W$, which can significantly affect performance, see \cite{EXTRA} for details.
	 
	 \section{Design of Algorithm}
	 \subsection{Proximal-Correction Algorithm}
	 In this section, we will construct an algorithm to solve \eqref{P}, which is inspired by a distributed optimization algorithm named EXTRA \cite{EXTRA}. 
	 EXTRA is to solve the following unconstrained distributed convex optimization problem:
	 \begin{equation*}
	 	\min_{z\in\mR^n} f(z)\triangleq\sum_{i=1}^Nf_i(z).
	 \end{equation*}
	 We introduce a collective variable $Z$ and write the gradient of $f$ at $Z$ as
	 \begin{align*}
	 	\nabla f(Z)=\left[
	 	\begin{array}{ccc}
	 		-& \left(\nabla f_1(z_1)\right)^\top &-\\
	 		-& \left(\nabla f_2(z_2)\right)^\top &-\\
	 		~& \vdots &~\\
	 		-& \left(\nabla f_N(z_N)\right)^\top &-
	 	\end{array}\right]\in\mR^{N\times n},\quad
 	    Z=\left[
 	    \begin{array}{ccc}
 	    	-& z_1^\top &-\\
 	    	-& z_2^\top &-\\
 	    	~& \vdots & ~\\
 	    	-& z_N^\top&-
 	    \end{array}\right]\in\mR^{N\times n},
	 \end{align*}
	 where $z_i\in\mR^n$ is a local copy of the global variable $z$, $i=1, 2,\dots, N$. 
	 Then the scheme of EXTRA is 
	 \begin{equation}\label{eq_EXTRA}
	 	\underbrace{Z^{k+1}=\tilde{W}Z^k-\alpha\nabla f(Z^k)}_{DGD}+
	 	\underbrace{\sum_{t=0}^{k}\left(W-\tilde{W}\right)Z^t}_{Correction},
	 	\quad \forall k\geq 0.
	 \end{equation}
	 As revealed in \cite{EXTRA}, this method can be viewed as a cumulative correction to DGD \cite{Nedic09b}, and then achieved a fast linear convergence rate. We extract this key idea of taking a cumulative correction and combine it with the proximal point algorithm (PPA) to develop our algorithm. Let us unfold the details.
	 
	 To find a solution to problem \eqref{P} in a centralized setting, Rockafellar \cite {Rockafellar76} proposed the proximal point algorithm (PPA):
	 \begin{equation*}
	 	z^{k+1}=\left(I+\alpha\tilde{T}\right)^{-1}(z^k), \quad\alpha>0,
	 \end{equation*}
     which is equivalent to
	 \begin{equation*}
	 	z^{k+1}=z^k-\alpha v^{k+1},\quad v^{k+1}\in\tilde{T}(z^{k+1}).
	 \end{equation*}
	 It is worth to note that the convergence of PPA is guaranteed for any value of the penalty parameter $\alpha$, as long as the operator $\tilde{T}$ is maximal monotone. For example, if $\tilde{T}$ is the subdifferential $\partial f$ of a lower semicontinuous convex function $f$, then the maximal monotonicity of $\tilde{T}$ is derived from the lower semicontinuity and convexity of $f$, thus the convergence of PPA is independent of the smoothness of the objective function, which differs from the standard gradient method.
	 
	 The authors of \cite{LiXX19} and \cite{LiXX21} employed PPA to construct the following distributed proximal point algorithm (abbreviated as DPPA):
	 \begin{equation*}
	 	z_i^{k+1}=\left(I+\alpha_k T_i\right)^{-1}\left(\sum_{j=1}^N\tilde{w}_{ij}z_j^k\right),
	 \end{equation*}
     or equivalently,
	 \begin{equation*}
	 	z_i^{k+1}=\sum_{j=1}^N\tilde{w}_{ij}z_j^k-\alpha_k v_i^{k+1},\quad v_i^{k+1}\in T_i(z_i^{k+1}).
	 \end{equation*}
	 In fact, \cite{LiXX19} and \cite{LiXX21} focused on the special cases $T_i=\partial f_i$ and $T_i=\partial L_i$, respectively, where $L_i$ is defined by \eqref{eq_Lag_Fun}. Here we provide a unified form of these two algorithms by using maximal monotone operators.
	 Let operator $T:\mR^{N\times n}\rightrightarrows\mR^{N\times n}$ satisfy
	 \begin{align*}
	 	T(Z)=\left[
	 	\begin{array}{ccc}
	 		-& \left(T_1(z_1)\right)^\top &-\\
	 		~& \vdots &~\\
	 		-& \left(T_N(z_N)\right)^\top &-
	 	\end{array}
	 	\right],\quad
	 	V=\left[
	 	\begin{array}{ccc}
	 		-& \left(v_1\right)^\top &-\\
	 		~& \vdots &~\\
	 		-& \left(v_N\right)^\top &-
	 	\end{array}
	 	\right]\in T(Z).
	 \end{align*}
	 Note that $T$ is the cartesian product of $N$ maximal monotone operators, and is also maximal monotone. Then DPPA can be rewritten into the following compact matrix form:
	 \begin{equation*}
	 	Z^{k+1}=\tilde{W}Z^k-\alpha_k V^{k+1},\quad V^{k+1}\in T(Z^{k+1}).
	 \end{equation*}
	  The work in \cite{LiXX19} and \cite{LiXX21} creatively extended PPA to distributed optimization problems. However, these two algorithms suggest adopting diminishing penalty parameters, i.e. $\alpha_k\to0$, $\sum_{k=0}^\infty\alpha_k=\infty$, to guarantee the convergence to an optimal solution, which leads to a slow convergence rate. This motivates us to accelerate these two algorithms. 
	  
	  Inspired by the fact that EXTRA introduced a cumulative correction for DGD and led to a significantly faster convergence rate, we attempt to add the same cumulative correction to DPPA and then propose the following distributed Proximal-Correction algorithm:
	 \begin{equation}\label{eq_alg}
	 	\underbrace{Z^{k+1}=\tilde{W}Z^k-\alpha V^{k+1}}_{DPPA}
	 	+\underbrace{\sum_{t=0}^k\left(W-\tilde{W}\right)Z^t}_{Correction},\quad \forall k\geq 0,
	 \end{equation}
	 where $V^k\in T(Z^k)$ for all $k\ge 1$. From the above formulation \eqref{eq_alg}, updating the decision vector $Z^{k+1}$ requires having access to all decision vectors from the initial step up to the $k$th step,
	 which is not convenient for practical calculation. We subtract the $k+1$ iteration from the $k+2$ iteration of Eq.\eqref{eq_alg} to yield 
	 \begin{subequations}\label{eq_compt}
	 	\begin{align}
	 		Z^1&=WX^0-\alpha V^1,\label{eq_compta}\\
	 		Z^{k+2}&=\left(I+W\right)Z^{k+1}-\tilde{W}Z^k-\alpha\left(V^{k+2}-V^{k+1}\right),\quad k\ge 0,
	 		\label{eq_comptb}
	 	\end{align}
	 \end{subequations}
	 Clearly, \eqref{eq_compt} cannot be implemented directly as it is, because we would need to know $Z^{k+2}$ to compute $V^{k+2}$, thus we introduce the following proximal mappings (resolvents) of the maximal monotone operators $T_i$ and $T$
	 \begin{equation*}
	 	\proxTi=\left(I+\alpha T_i\right)^{-1},\quad\proxT=\left(I+\alpha T\right)^{-1},\quad\alpha>0,
	 \end{equation*}
	 where $I$ represents the identity operator. Note that the proximal mappings are single-valued and non-expansive. From \eqref{eq_compt}, we have 
	 \begin{align*}
	 	&WZ^0\in(I+\alpha T)(Z^1),\\
	 	(I+W)Z^{k+1}-&\tilde{W}Z^k+\alpha V^{k+1}\in\left(I+\alpha T\right)\left(Z^{k+2}\right),
	 \end{align*}
	 which gives the update formulas of Proximal-Correction as
	 \begin{subequations}\label{eq_itera_PC}
	 \begin{align}
	 	Z^1&=\proxT\left(WZ^0\right), \label{eq_itera_PCa}\\
	 	V^1&=\left(WZ^0-Z^1\right)/\alpha,\label{eq_itera_PCb}\\
	  Z^{k+2}&=\proxT\left((I+W)Z^{k+1}-\tilde{W}Z^k+\alpha V^{k+1}\right),\label{eq_itera_PCc}\\
	  V^{k+2}&=\left(\left(I+W\right)Z^{k+1}-\tilde{W}Z^k-Z^{k+2}+\alpha V^{k+1}\right)\big/\alpha,\quad\text{for}\quad k=0,1,2,\dots\label{eq_itera_PCd}
	  \end{align}
	 \end{subequations}
	 It is known that if $T_i=\partial f_i$ with a convex function $f_i$, then 
	 \begin{equation}\label{exam_1}
	 	u=\proxTi(z)\Longleftrightarrow 
	 	u=\mathop{\arg\min}\limits_{\zeta}\left\{f_i(z)+\frac{1}{2\alpha}\Vert\zeta-z\Vert^2\right\},
	 \end{equation}
	 and if $T_i=\partial L_i$ with a convex-concave function $L_i(\cdot,\cdot)$, then 
	 \begin{equation}\label{exam_2}
	 	(u,v)=\proxTi(x,y)\Leftrightarrow (u,v)=\mathop{\arg\min\max}\limits_{\zeta, \eta}
	 	\left\{L_i(x,y)+\frac{1}{2\alpha}\Vert\zeta-x\Vert^2-\frac{1}{2\alpha}\Vert\eta-y\Vert^2\right\}.
	 \end{equation}
	 Now we give the implementation steps of the distributed Proximal-Correction algorithm as follows:
    \begin{algorithm}[htbp]
    	\caption{Distributed Proximal-Correction Algorithm}
    	\label{alg_1}
    	Choose penalty parameter $\alpha>0$, and mixing matrices $W, \tilde{W}\in\mR^{n\times n}$;\\
    	Pick any $z_i^0\in\mR^n, i=1, 2,\dots, N$;
    	\begin{algorithmic}
    		\STATE $z_i^1=\proxTi\left(\sum_{i=1}^Nw_{ij}z_j^0\right)$;\\
    		$v_i^1=\left(\sum_{i=1}^Nw_{ij}z_j^0-z_i^1\right)\big/\alpha$;
    		\FOR{$k=0,1,2,\dots$} 
    		\vspace{-1em}
    		\item {
    			\begin{align*}
    				z_i^{k+2}&=\proxTi\left(z_i^{k+1}+\sum_{j=1}^Nw_{ij}z_j^{k+1}-
    				\sum_{j=1}^N\tilde{w}_{ij}z_j^k+\alpha v_i^{k+1}\right);\hspace{15em}\\
    				v_i^{k+2}&=\left(z_i^{k+1}+\sum_{j=1}^Nw_{ij}z_j^{k+1}-
    				\sum_{j=1}^N\tilde{w}_{ij}z_j^k-z_i^{k+2}+\alpha v_i^{k+1}\right)\bigg/\alpha;
    			\end{align*}
    		}
    		\ENDFOR
    	\end{algorithmic}
    \end{algorithm}
    \subsection{Inexact Computation of the Proximal Step}
    	From a practical point of view the computation of the proximal mapping $\proxTi$ can be as difficult as the computation of a solution to its equivalent optimization problem (e.g., see \eqref{exam_1} and \eqref{exam_2}), even though the strong convexity of the latter can be a key advantage, so it is essential to replace the equations 
    \begin{align*}
    	Z^1&=\proxT\left(WZ^0\right), \\
    	Z^{k+2}&=\proxT\left(\left(I+W\right)Z^{k+1}-\tilde{W}Z^k+\alpha V^{k+1}\right)
    \end{align*}
    by looser relations 
    \begin{align*}
    	Z^1&\approx\proxT\left(WZ^0\right), \\
    	Z^{k+2}&\approx\proxT\left(\left(I+W\right)Z^{k+1}-\tilde{W}Z^k+\alpha V^{k+1}\right), 
    \end{align*}
    respectively. 
    In this paper, we consider two general criteria for the approximate calculation of the proximal mapping such that the looser relations are computationally implementable for a wide variety of problems. These two criteria are first presented by Rockafellar in \cite{Rockafellar76}, and they are
    	\begin{align*}
    	&({\rm A1})\hspace{5em}
    	\left\Vert Z^{k+2}-\proxT\left(\widehat{Z}^{k+1}\right)\right\Vert\le\varepsilon_{k+1}, \quad
    	\sum_{k=0}^\infty\varepsilon_{k+1}<\infty,\\
    	&({\rm B1})\hspace{5em}
    	\left\Vert Z^{k+2}-\proxT\left(\widehat{Z}^{k+1}\right)\right\Vert\le\delta_{k+1}\Vert Z^{k+2}-Z^{k+1}\Vert, \quad 
    	\sum_{k=0}^\infty\delta_{k+1}<\infty,
    \end{align*}
    where 
    \begin{equation}\label{eq_zhat}
    	\widehat{Z}^{k+1}=\left(I+W\right)Z^{k+1}-\tilde{W}Z^k+\alpha V^{k+1}, \quad\forall k\ge0.
    \end{equation}
    In practical calculations, (A1) (or (B1)) is implemented by agents in a distributed manner, here we consider the compact form for the convenience of analysis.
    For the approximate calculation of $\proxT\left(\widehat{Z}^{k+1}\right)$, we have the estimate (Proposition 3 \cite{Rockafellar76})
    \begin{equation}\label{ine_estm}
    	\left\Vert Z^{k+2}-\proxT\left(\widehat{Z}^{k+1}\right)\right\Vert\le\alpha{\rm dist}
    	\left(\0,S_T^{k+1}\left(Z^{k+2}\right)\right),
    \end{equation} 
    where 
    \begin{equation}\label{eq_S_T}
    	S_T^{k+1}\left(Z\right)=T(Z)+\left(Z-\widehat{Z}^{k+1}\right)\big/\alpha.
    \end{equation}
    Therefore, (A1) and (B1) are respectively implied by the following criteria 
    \begin{align*}
    	&({\rm A2})\hspace{5em}
    	{\rm dist}\left(\0,S_T^k(Z^{k+2})\right)\le\varepsilon_{k+1}/\alpha,\quad
    	\sum_{k=0}^\infty\varepsilon_{k+1}<\infty;\\
    	&({\rm B2})\hspace{5em}
    	{\rm dist}\left(\0,S_T^k(Z^{k+2})\right)\le\left(\delta_{k+1}/\alpha\right)\Vert Z^{k+2}-Z^{k+1}\Vert,\quad
    	\sum_{k=0}^\infty\delta_{k+1}<\infty.
    \end{align*}
    Note that each listed criterion with $\varepsilon_{k+1}=0$ or $\delta_{k+1}=0$ implies the exact case $$Z^{k+2}=\proxT\left(\widehat{Z}^{k+1}\right),$$ 
    therefore the theoretical results that will be shown in Section 4 are valid for Algorithm \ref{alg_1}.
    Criterion $({\rm A2})$ can be rewritten in another way. In fact,
    \begin{align}
    	&{\rm dist}\left(\0,S_T(Z^{k+2})\right)\le\varepsilon_{k+1}/\alpha\notag\\
    	\Longleftrightarrow& \ \exists \ e^{k+1}, \ \Vert e^{k+1}\Vert\le\varepsilon_{k+1}/\alpha, \ e^{k+1}\in S_T(Z^{k+2})\notag\\
    	\Longleftrightarrow& \ \exists \ e^{k+1}, \ \Vert e^{k+1}\Vert\le\varepsilon_{k+1}, \ \left(\widehat{Z}^{k+1}-Z^{k+2}+e^{k+1}\right)\big/\alpha\in T(Z^{k+2})\notag\\
    	\Longleftrightarrow& \ \exists \ e^{k+1}, \ \Vert e^{k+1}\Vert\le\varepsilon_{k+1}, \ Z^{k+2}=\proxT\left(\widehat{Z}^{k+1}+e^{k+1}\right).\label{eq_A2}
    \end{align}
    
    \begin{rem}\label{rem_1}
    Criterion $({\rm A2})$ is generally convenient as a measure of how near $Z^\prime$ approximates $\proxT\left(Z\right)$. 
    Take problem \eqref{DP} as an example, we have $T_i=\partial L_i$, where $L_i$ is defined by \eqref{eq_Lag_Fun}. Suppose that functions $f_i$ and $h_i$ are differentiable. Denote
    \begin{equation*}
    	\tilde{L}_i(x,\lambda)=f_i(x)+\lambda^\top h_i(x), \quad \tilde{\Omega}_i=\Omega_i\times K^\circ,
    \end{equation*}
     then we have 
     \begin{equation*}
     	\partial L_i(x,\lambda)=\nabla\tilde{L}_i(x,\lambda)+N_{\tilde{\Omega}_i}(x,\lambda),
     \end{equation*}
    where $N_{\tilde{\Omega}_i}(x,\lambda)$ represents the normal cone of the subset $\tilde{\Omega}_i$ at $(x,\lambda)$.
    In this case, criterion $({\rm A2})$ for $(y,\mu)\approx\proxTi(x,\lambda)$ is 
    \begin{align*}
     &{\rm dist}\left(\0,S_{T_i}(y,\mu)\right)\le\varepsilon/\alpha\\
     \Longleftrightarrow \ &{\rm dist}\left(\0,\nabla\tilde{L}_i(y,\mu)+N_{\tilde{\Omega}_i}(y,u)+[(y,\mu)-(x,\lambda)]/\alpha\right)\le\varepsilon/\alpha\\
     \Longleftrightarrow \ &{\rm dist}\left(-\nabla\tilde{L}_i(y,\mu)-[(y,\mu)-(x,\lambda)]/\alpha,
     N_{\tilde{\Omega}_i}(y,\mu)\right)\le\varepsilon/\alpha\\
     \Longleftrightarrow \ & \left\Vert w-{\rm Proj}_{N_{\tilde{\Omega}_i}(y,\mu)}(w)\right\Vert\le\varepsilon/\alpha,
    \end{align*} 
    where $w= -\nabla\tilde{L}_i(y,\mu)-[(y,\mu)-(x,\lambda)]/\alpha$, and $S_{T_i}(y,\mu)=T_i(y,\mu)+[(y,\mu)-(x,\lambda)]/\alpha$.  While criterion $(\rm A1)$ cannot be directly verified to be satisfied when we don't have the exact value of $\proxTi(x,\lambda)$.
    \end{rem}
    In these criteria, the error sequence $\{\varepsilon_k\}_{k\ge0}$ (or $\{\delta_k\}_{k\ge0}$) is required to be summable. If this condition is relaxed to only $\varepsilon_k\to 0$ (as $k\to\infty$), convergence to a solution to problem \eqref{P} may fail. A counterexample will be shown in numerical simulations.
    
    \section{Convergence Analysis}
    In this section, we will analyze the convergence (Theorem \ref{thm_2}) of Proximal-Correction Algorithm at first. Furthermore, when the proximal steps of Proximal-correction Algorithm are inexactly calculated, we will show that criterion (A1) (or (A2)) is enough to guarantee the convergence (Theorem \ref{thm_3}) to a solution to problem \eqref{P}, even though there may be more than one solution; when using criterion (B1) (or (B2)) to compute the proximal steps, under the assumption of $\Phi^{-1}$, the inverse of a
    constructed operator $\Phi$ defined by \eqref{eq_phi}, being Lipschitz continuous at $\0$, it will be shown that the convergence rate is linear (Theorem \ref{thm_4}).
    
    Under Assumption \ref{ass_2}, the matrix $$C=\tilde{W}-W$$ is positive semidefinite. Introducing an auxiliary sequence
    \begin{equation}\label{eq_Y}
    	Y^k\triangleq\sum_{t=0}^k\sqrt{C}Z^t,\quad \forall k\ge0,
    \end{equation} 
     then \eqref{eq_alg} can be equivalently written as
     \begin{subequations}\label{eq_PD}
     	\begin{align}
     		Z^{k+1}&=\tilde{W}Z^k-\alpha V^{k+1}-\sqrt{C}Y^k,\label{eq_PDa}\\
     		Y^{k+1}&=Y^k+\sqrt{C}Z^{k+1},\label{eq_PDb}
     	\end{align}
     \end{subequations}
    where $V^{k+1}\in T(Z^{k+1})$ for all $k\ge0$.
    We first give a lemma that motivates us to perform the above equivalent transformation for Proximal-Correction. 
    \begin{lem}\label{lem_1}
    	 Under Assumption \ref{ass_2}, for any fixed scalar $\alpha>0$, if $z^*$ is a solution to problem \eqref{P}, then there exist $Y^*=\sqrt{C}Y$ for some $Y\in\mR^{N\times n}$ and $Z^*={\bf 1}\otimes \left(z^*\right)^\top$ such that
    	\begin{subequations}\label{eq_lem1}
    		\begin{align}
    			\0&\in \alpha T(Z^*)+\sqrt{C}Y^*,\label{eq_lem1a}\\
    			\0&=-\sqrt{C}Z^*.\label{eq_lem1b}
    		\end{align}
    	\end{subequations}
    Conversely, if there exists a pair $(Z^*,Y^*)\in\mR^{N\times n}\times\mR^{N\times n}$ satisfying \eqref{eq_lem1}, then $Z^*=\1\otimes \left(z^*\right)^\top$, and $z^*$ is a solution to problem \eqref{P}.
    \end{lem}
    {\rm\bf Proof}  According to Assumption \ref{ass_2} and the definition of $C$, we have
    \begin{equation*}
    	{\rm null}\{\sqrt{C}\}={\rm null}\{\tilde{W}-W\}={\rm span}\{\1\}.
    \end{equation*}
    Hence $Z^*\in\mR^{N\times n}$ is consensual if and only if \eqref{eq_lem1b} holds.

    Suppose $z^*$ is a solution to problem \eqref{P}, i.e. $0\in\sum_{i=1}^NT_i(z^*)$. Take
    $Z^*=\1\otimes \left(z^*\right)^\top$, then \eqref{eq_lem1b} is satisfied, and we have
    $\0\in\1^\top T(Z^*)$, i.e. $\1^\top V^*=\0$ for some $V^*\in T(Z^*)$. Since ${\rm span}\left\{\sqrt{C}\right\}=\left({\rm span}\{\1\}\right)^\perp$, 
    we have $V^*\in{\rm span}\left\{\sqrt{C}\right\}$, which implies that there exists $Y\in\mR^{N\times n}$ with $\alpha V^*=-\sqrt{C}Y$. Let $Y^*={\rm Proj}_{\sqrt{C}}\left(Y\right)$, then $\sqrt{C}Y^*=\sqrt{C}Y$, and \eqref{eq_lem1a} holds.
    
    Conversely, suppose that there is a pair $(Z^*,Y^*)$ such that \eqref{eq_lem1} holds, then $Z^*=\1\otimes \left(x^*\right)^\top$ is directly derived from \eqref{eq_lem1b} for some $z^*\in\mR^n$. Under Assumption \ref{ass_2}(b), we have $\1^\top C=\0$ and 
    $\1^\top\sqrt{C}=\0$. Multiplying \eqref{eq_lem1a} by $\1^\top$ gives $\0\in\1^\top T(Z^*)$, i.e.,
    $\0\in\sum_{i=1}^NT_i(z^*)$.
    \qed\\

    From Lemma \ref{lem_1}, solving problem \eqref{P} turns into finding a pair $(Z^*,Y^*)$ to meet \eqref{eq_lem1}. We introduce an operator  
    $\Phi:\mR^{N\times n}\times\mR^{N\times n}\rightrightarrows\mR^{N\times n}\times\mR^{N\times n}$ satisfying 
    \begin{equation}\label{eq_phi}
    \Phi(Z,Y)=\left[
    	\begin{array}{c}
    		\alpha T(Z)+\sqrt{C}Y\\
    		-\sqrt{C}Z
    	\end{array}\right].
    \end{equation}
    It can be asserted that the operator $\Phi$ is maximal monotone, we demonstrate the details with the following Minty theorem.
    \begin{thm}[Minty Theorem \cite{Rockafellar09}]\label{thm_1}
    	A monotone operator $\Phi$ is maximal if and only if the operator $I+\Phi$ is surjective.
    \end{thm}
    \begin{propo}\label{propo_1}
    	The operator $\Phi$ is maximal monotone.
    \end{propo}
    {\rm \bf Proof} We first verify the monotonicity of $\Phi$. For any
    $\omega_1\in\Phi(Z_1,Y_1)$, $\omega_2\in\Phi(Z_2,Y_2)$, we have 
    \begin{align*}
    	\omega_1=\left[
    	\begin{array}{c}
    		\alpha V_1+\sqrt{C}Y_1\\
    		-\sqrt{C}Z_1
    	\end{array}\right],\quad
    \omega_2=\left[
    \begin{array}{c}
    	\alpha V_2+\sqrt{C}Y_2\\
    	-\sqrt{C}Z_2
    \end{array}\right],\quad V_i\in T(Z_i), i=1, 2.
    \end{align*}
   Then 
   \begin{align*}
   	&\left\langle\omega_1-\omega_2, (Z_1,Y_1)-(Z_2,Y_2)\right\rangle\\
   	=&\langle\alpha(V_1-V_2)+\sqrt{C}(Y_1-Y_2), Z_1-Z_2\rangle+\langle-\sqrt{C}(Z_1-Z_2), Y_1-Y_2\rangle\\
   	=&\alpha\langle V_1-V_2, Z_1-Z_2\rangle\\
   \ge&0.
   \end{align*}
    Therefore, the monotonicity of $\Phi$ follows from $T$. 
    
     As for the maximality, according to Theorem \ref{thm_1}, we need to verify the operator $I+\Phi$ is surjective. Under Assumption \ref{ass_1}, for any $\alpha>0$, $\alpha T$ is maximal monotone, it follows that $\alpha(I+C)^{-1}T$ is also maximal monotone with  $I+C\succ 0$. Then the operator $I+\alpha(I+C)^{-1}T$ is surjective. Take any $(\zeta,\eta)\in\mR^{N\times n}\times\mR^{N\times n}$, suppose that there is a pair $(Z,Y)$ such that 
    $(\zeta,\eta)\in(I+\alpha\Phi)(Z,Y)$, i.e.,
    \begin{subequations}
    	\begin{align}
    		&\zeta\in Z+\alpha T(Z)+\sqrt{C}Y,\label{eq_lem2a}\\
    		&\eta=Y-\sqrt{C}Z.\label{eq_lem2b}
    	\end{align}
    \end{subequations}
    Eliminating the variable $Y$ from \eqref{eq_lem2a} yields
    \begin{equation}\label{eq_in}
    	(I+C)^{-1}\left(\zeta-\sqrt{C}\eta\right)\in\left(I+\alpha(I+C)^{-1}T\right)(Z).
    \end{equation}
    In view of that $I+\alpha(I+C)^{-1}T$ is surjective, there does exist $Z$ satisfying \eqref{eq_in}, i.e.,
    \begin{align*}
    	Z&=\left(I+\alpha(I+C)^{-1}T\right)^{-1}\left((I+C)^{-1}(\zeta-\sqrt{C}\eta)\right),\\
    	Y&=\eta+\sqrt{C}Z.
    \end{align*}
    Therefore, the operator $I+\Phi$ is surjective, and the proof is completed.\qed\\

    Now we present the first main theoretical result, which establishes the convergence of the proposed Proximal-Correction algorithm to a solution to problem \eqref{P}.
    \subsection*{Convergence of Proximal-Correction Algorithm}
    \begin{thm}\label{thm_2}
    	Under Assumptions \ref{ass_1} and \ref{ass_2}, for any $\alpha>0$, all sequences $\{z_i^k\}_{k\ge0}, i=1, 2,\dots, N$, generated by Algorithm \ref{alg_1} converge to the same solution $z^\infty$ to problem \eqref{P}, and for each $i=1, 2,,\dots, N$, the sequence $\{v_i^k\}_{k\ge0}$ converges to an element of $T_i(z^\infty)$.
    \end{thm}
    {\rm \bf Proof} \ We deal with the equivalent transformation \eqref{eq_PD} of the iterative scheme \eqref{eq_itera_PC} of Proximal-Correction Algorithm. The idea of the proof is as follows: (1) We would first show that there exists a cluster point $(Z^\infty,Y^\infty)$ of $\{(Z^k,Y^k)\}_{k\ge0}$ satisfying $\0\in\Phi(Z^\infty,Y^\infty)$; (2) then we argue that $\{\left(Z^k,Y^k\right)\}_{k\ge0}$ cannot have more than one cluster; (3) The claimed results are obtained by applying Lemma \ref{lem_1}.
    
    {\bf Step 1:} From Eq.\eqref{eq_PD}, we have
    \begin{align*}
    	\left[
    	\begin{array}{c}
    		Z^{k+1}\\
    		Y^{k+1}
    	\end{array}\right]
        &=\left[
        \begin{array}{c}
        	\tilde{W}Z^k\\
        	Y^k
        \end{array}\right]
       -\left[
       \begin{array}{c}
       	\alpha V^{k+1}+\sqrt{C}Y^k\\
       	-\sqrt{C}Z^{k+1}
       \end{array}\right]\\
       &=\left[
       \begin{array}{c}
       	\tilde{W}Z^k\\
       	Y^k
       \end{array}\right]
       -\left[
       \begin{array}{c}
       	\alpha V^{k+1}+\sqrt{C}Y^{k+1}\\
       	-\sqrt{C}Z^{k+1}
       \end{array}\right]
       +\left[
       \begin{array}{c}
       	\sqrt{C}(Y^{k+1}-Y^k)\\
       	\0
       \end{array}\right]\\
       &=\left[
       \begin{array}{c}
       	\tilde{W}Z^k\\
       	Y^k
       \end{array}\right]
       -\left[
       \begin{array}{c}
       	\alpha V^{k+1}+\sqrt{C}Y^{k+1}\\
       	-\sqrt{C}Z^{k+1}
       \end{array}\right]
       +\left[
       \begin{array}{c}
       	CZ^{k+1}\\
       	\0
       \end{array}\right]
    \end{align*}
    which gives
    \begin{align*}
    	\left[
    	\begin{array}{cc}
    		I-C &\0\\
    		\0  &I
    	\end{array}\right]
    \left[
    \begin{array}{c}
    	Z^{k+1}\\
    	Y^{k+1}
    \end{array}\right]
    =\left[
    \begin{array}{cc}
    	\tilde{W} &\0\\
    	\0  &I
    \end{array}\right]
    \left[
    \begin{array}{c}
    	Z^k\\
    	Y^k
    \end{array}\right]
    -\left[
    \begin{array}{c}
    	\alpha V^{k+1}+\sqrt{C}Y^{k+1}\\
    	-\sqrt{C}Z^{k+1}
    \end{array}\right].
    \end{align*}
    For convenience, denote 
    \begin{equation}\label{eq_def_pq}
    \begin{aligned}
    	&P=\left[
    	\begin{array}{cc}
    		I-C &\0\\
    		\0  &I
    	\end{array}\right]\succ\0,\quad
        Q=\left[
        \begin{array}{cc}
        	\tilde{W} &\0\\
        	\0  &I
        \end{array}\right]\succ\0,\\
    &\xi^k=\left[
    \begin{array}{c}
    	Z^k\\
    	Y^k
    \end{array}\right],\quad
    \omega^{k+1}=\left[
    \begin{array}{c}
    	\alpha V^{k+1}+\sqrt{C}Y^{k+1}\\
    	-\sqrt{C}Z^{k+1}
    \end{array}\right]\in\Phi(\xi^{k+1}).
    \end{aligned}
    \end{equation}
    The positive definite properties of $P$ and $Q$ are guaranteed by Assumption \ref{ass_2}. Then we obtain
    \begin{equation}\label{eq_key}
    	P\xi^{k+1}+\omega^{k+1}=Q\xi^k, \quad \omega^{k+1}\in\Phi(\xi^{k+1}).
    \end{equation}
    Let $z^*$ be an arbitrary solution to problem \eqref{P}, and $Z^*=\1\otimes\left(x^*\right)^\top$. From Lemma \ref{lem_1}, there exists $Y^*\in\mR^{N\times n}$ such that $0\in\Phi(Z^*,Y^*)$. Denote
    $\xi^*=\left((Z^*)^\top,(Y^*)^\top\right)^\top$, due to the doubly-stochasticity of the mixing matrices $W$ and $\tilde{W}$, it holds
    \begin{equation*}
    	P\xi^*=\xi^*,\quad Q\xi^*=\xi^*, \quad\text{and}\quad\0\in\Phi(\xi^*).
    \end{equation*}
    Subtracting $\xi^*$ at both sides of \eqref{eq_key} yields
    \begin{equation*}
    	P\left(\xi^{k+1}-\xi^*\right)+\omega^{k+1}=Q\left(\xi^k-\xi^*\right),
    \end{equation*}
    and then taking the squared norm to give
    \begin{equation*}
    	\left\Vert\xi^{k+1}-\xi^*\right\Vert^2_{P^2}+2\langle P(\xi^{k+1}-\xi^*),\omega^{k+1}-\0\rangle
    	+\Vert\omega^{k+1}\Vert^2=\left\Vert\xi^k-\xi^*\right\Vert^2_{Q^2}.
    \end{equation*}
    The monotonicity of $\Phi$ implies 
    \begin{equation*}
    	\langle \xi^{k+1}-\xi^*,\omega^{k+1}-\0\rangle\ge0.
    \end{equation*}
    Furthermore, considering the matrix $P\succ0$, it holds that 
    \begin{equation*}
    	\langle P(\xi^{k+1}-\xi^*),\omega^{k+1}-\0\rangle\ge0.
    \end{equation*} 
     Thus we have
    \begin{equation}\label{ieq_cruc1}
    	\left\Vert\xi^{k+1}-\xi^*\right\Vert^2_{P^2}+\Vert\omega^{k+1}\Vert^2\le
    	\left\Vert\xi^k-\xi^*\right\Vert^2_{Q^2}.
    \end{equation}
   Note that
    \begin{equation*}
    	P^2-Q^2=\left[
    	\begin{array}{cc}
    		(I+W)(I+W-2\tilde{W}) &\0\\
    		\0 & \0
    	\end{array}\right],
    \end{equation*}
    Under Assumption \ref{ass_2}, it is clear that $P^2-Q^2\succcurlyeq 0$, which follows 
    \begin{equation}\label{ieq_cruc2}
    	\left\Vert\xi^{k+1}-\xi^*\right\Vert^2_{P^2}+\Vert\omega^{k+1}\Vert^2\le
    	\left\Vert\xi^k-\xi^*\right\Vert^2_{P^2}.
    \end{equation}
    Summing \eqref{ieq_cruc2} over $k$ from $0$ to $\infty$ yields
    \begin{equation}\label{ieq_omega}
    	\sum_{k=0}^\infty\Vert\omega^{k+1}\Vert^2\le \Vert\xi^0-\xi^*\Vert^2_{P^2}<\infty.
    \end{equation}
    The sequence $\left\{\Vert\xi^k-\xi^*\Vert^2_{P^2}\right\}_{k\ge0}$ is decreasing and bounded  below, hence the limit 
    \begin{equation*}
    	\lim_{k\to\infty}\Vert\xi^k-\xi^*\Vert_{P^2}
    \end{equation*}
    exists and the sequence $\{\xi^k\}_{k\ge0}$ is bounded. From \eqref{ieq_omega} we have
    \begin{equation*}
    	\lim_{k\to\infty}\Vert\omega^k\Vert=0.
    \end{equation*}
    Let $\xi^\infty$ be an arbitrary cluster point of $\{\xi^k\}_{k\ge0}$, for any
    $\omega\in\Phi(\xi)$, the monotonicity of $\Phi$ gives
    \begin{equation*}
    	\langle\omega-\omega^{k+1}, \xi-\xi^{k+1}\rangle\ge0.
    \end{equation*}
    Taking $k\to\infty$ yields
    \begin{equation*}
    	\langle\omega-\0, \xi-\xi^\infty\rangle\ge0 \quad \text{ for any}\quad \xi, \omega\quad\text{satisfying}\quad \omega\in\Phi(\xi),
    \end{equation*}
    which implies $\0\in\Phi(\xi^\infty)$, in view of the maximality of $\Phi$.\\
    
    {\bf Step 2:} This step is to show that there cannot be more than one cluster point of $\{\xi^k\}_{k\ge0}$. Suppose
    there are two cluster points: $\xi^\infty_1\neq\xi^\infty_2$. Then $\0\in\Phi(\xi^\infty_i)$ for $i=1, 2$, as just seen, so that each $\xi^\infty_i$ can play the role of $\xi^*$ in \eqref{ieq_cruc2}, and  we get the existence of the limits
    \begin{equation}\label{eq_conv}
    	\lim_{k\to\infty}\Vert\xi^k-\xi^\infty_i\Vert^2_{P^2}=\mu_i<\infty,\quad \text{for}\quad i=1, 2.
    \end{equation}
    Writing 
    \begin{equation*}
    	\Vert\xi^k-\xi^\infty_2\Vert^2_{P^2}=\Vert\xi^k-\xi^\infty_1\Vert^2_{P^2}+
    	2\langle P(\xi^k-\xi^\infty_1),P(\xi^\infty_1-\xi^\infty_2)\rangle
    	+\Vert\xi^\infty_1-\xi^\infty_2\Vert^2_{P^2},
    \end{equation*}
    we see that the limit of $2\langle P(\xi^k-\xi^\infty_1),P(\xi^\infty_1-\xi^\infty_2)\rangle$ must also 
    exist and
    \begin{equation*}
    	2\lim_{k\to\infty}\langle P(\xi^k-\xi^\infty_1),P(\xi^\infty_1-\xi^\infty_2)\rangle
    	=\mu_2-\mu_1-\Vert\xi^\infty_1-\xi^\infty_2\Vert^2_{P^2} .
    \end{equation*}
    But this limit cannot be different from $0$, because $\xi^\infty_1$ is a cluster point of
    $\{\xi^k\}_{k\ge0}$. Therefore
    \begin{equation*}
    	\mu_2-\mu_1=\Vert\xi^\infty_1-\xi^\infty_2\Vert^2_{P^2}>0.
    \end{equation*}
    However, the same argument works with $\xi^\infty_1$ and $\xi^\infty_2$ reversed, so that also
    $\mu_1-\mu_2>0$. This is a contradiction which establishes the uniqueness of the cluster point 
    $\xi^\infty$. 
    
    {\bf Step 3:} Now it can be asserted that there exists
     $\xi^\infty=\left((Z^\infty)^\top,(Y^\infty)^\top\right)^\top\in\mR^{N\times n}\times\mR^{N\times n}$
     such that 
     \begin{subequations}
     	\begin{align}
     	&\lim_{k\to\infty}\Vert Z^k-Z^\infty\Vert=0,\quad\lim_{k\to\infty}\Vert Y^k-Y^\infty\Vert=0,\\
     	&\0\in\Phi(Z^\infty,Y^\infty)\label{eq_zopt}.
     	\end{align}
     \end{subequations}
     From Lemma \ref{lem_1}, we obtain that $Z^\infty=\1\otimes(z^\infty)^\top$ and $z^\infty$ is a solution to problem \eqref{P}.
     
     The rest is to show that $\{v_i^k\}_{k\ge0}$ converges to an element of $T_i(z^\infty)$ for each $i=1, 2,\dots, N$. Because of $\lim_{k\to\infty}\Vert\omega^k\Vert=0$, we have 
     \begin{equation*}
     	\lim_{k\to\infty}\left\Vert\alpha V^{k+1}+\sqrt{C}Y^{k+1}\right\Vert=0.
     \end{equation*}
    Since $Y^k\to Y^\infty$, it is clear that there exists $V^\infty\in\mR^{N\times n}$ such that
    \begin{equation}
    \lim_{k\to\infty}\Vert V^{k+1}-V^\infty\Vert=0 \quad\text{and}\quad \alpha V^\infty+\sqrt{C}Y^\infty=\0.
    \end{equation}
    where 
    \begin{align*}
    	V^\infty=\left[
    	\begin{array}{c}
    		(v^\infty_1)^\top\\
    		\vdots\\
    		(v^\infty_N)^\top
    	\end{array}\right]\in\mR^{N\times n}.
    \end{align*}
     Suppose that $V^\infty\notin T(X^\infty)$, i.e.
    \begin{align*}
    	\0=\left[
    	\begin{array}{c}
    		\alpha V^\infty+\sqrt{C}Y^\infty\\
    		-\sqrt{C}Z^\infty
    	\end{array}\right]\notin\Phi(Z^\infty,Y^\infty),
    \end{align*}
    which gives a contradiction to \eqref{eq_zopt}, thus it must hold $V^\infty\in T(Z^\infty)$, i.e., 
     \begin{equation*}
    	\lim_{k\to\infty}\Vert v_i^k-v^\infty_i\Vert=0 \quad\text{and} \quad v^\infty_i\in T_i(z^\infty),\quad
    	 \text{for each}\quad i=1, 2,\dots,N.
    \end{equation*}\qed
    \subsection*{Convergence of Inexact Proximal-Correction Algorithm}
     By introducing the auxiliary sequence $\{\xi^k\}_{k\ge0}$ and maximal operator $\Phi$, we obtain Eq.\eqref{eq_key}, an equivalent transformation of the scheme of Proximal-Correction. Eq.\eqref{eq_key} implies 
     \begin{equation}\label{eq_xi}
     	\xi^{k+1}=\left(I+P^{-1}\Phi\right)^{-1}\left(P^{-1}Q\xi^k\right)=\proxtP\left(\xi^k\right),
     \end{equation}  which is close to the scheme of the proximal point method for finding a point $\xi$ to meet $\0\in\Phi(\xi)$. 
    Naturally, when the proximal mapping $\proxT$ is  approximately computed using criteria (A1), (A2), (B1), and (B2), respectively, we would like to construct corresponding criteria for
    \begin{equation*}
    	\xi^{k+1}\approx\proxtP(\xi^k) 
    \end{equation*} 
    to develop the convergence analysis of the inexact Proximal-Correction Algorithm.
    
    From \eqref{eq_itera_PCb} and \eqref{eq_itera_PCd}, we have
    \begin{equation*}
    	Z^{k+2}-Z^{k+1}=WZ^{k+1}-\tilde{W}Z^k-\alpha\left(V^{k+2}-V^{k+1}\right),
    \end{equation*}
    summing over $k$ gives
    \begin{equation*}
    	Z^{k+2}=\tilde{W}Z^{k+1}-\alpha V^{k+2}+\sum_{k=0}^{k+1}\left(W-\tilde{W}\right)Z^t,\quad\forall k\ge0,
    \end{equation*}
    which means that \eqref{eq_alg} is still valid, but $V^{k+1}\notin T(Z^{k+1})$ for all $k\ge0$, because  \eqref{eq_itera_PCa} and \eqref{eq_itera_PCc} are no longer satisfied. 
    By the definitions of $Y^k$ and $\widehat{Z}^{k+1}$, we have
    \begin{align}
    	\widehat{Z}^{k+1}&=WZ^{k+1}+Z^{k+1}-\left(\tilde{W}Z^k-\alpha V^{k+1}\right)\notag\\
    	&=WZ^{k+1}+\sum_{t=0}^k\left(W-\tilde{W}\right)Z^t\notag\\
    	&=\tilde{W}Z^{k+1}+\sum_{t=0}^{k+1}\left(W-\tilde{W}\right)Z^t\notag\\
    	&=\tilde{W}Z^{k+1}-\sqrt{C}Y^{k+1}.\label{eq_zhat_y}
    \end{align}
    The sequence $\{V^k\}$ can be computed by using $\{Z^k\}$ and $\{Y^k\}$, i.e.
    \begin{align}\label{eq_VZY}
    	V^{k+1}&=\left(\tilde{W}Z^k+\sum_{t=0}^k\left(W-\tilde{W}\right)Z^t-Z^{k+1}\right)\bigg/\alpha\notag\\
    	&=\left(\tilde{W}Z^k-\sqrt{C}Y^k-Z^{k+1}\right)\big/\alpha,\quad \forall k\ge0. \hspace{4.5em}
    \end{align}
    Therefore, the iterative scheme of the inexact Proximal-Correction Algorithm can be reformulated as 
    \begin{subequations}\label{itera_yPC}
    	\begin{align}
    		Z^{k+1}&\approx\proxT\left(\tilde{W}Z^k-\sqrt{C}Y^k\right)\label{itera_yPCa}\\
    		Y^{k+1}&=Y^k+\sqrt{C}Z^{k+1},\label{itera_yPCb}
    	\end{align}
    \end{subequations}
    where the initial point $(Z^0,Y^0)$ satisfies $Y^0=\sqrt{C}Z^0$.
    To be clear, we would like to estimate the error of $\xi^{k+1}$ approximation to $\proxtP(\xi^k)$, which requires the following proposition.
    \begin{propo}\label{propo_2}
    	Suppose that 
    	\begin{equation*}
    		Z^{k+2}=\proxT\left(\widehat{Z}^{k+1}\right)+E^{k+1},
    	\end{equation*}
    then we have
    \begin{enumerate}[(a)]
    	\item {
         \begin{equation*}
         	\xi^{k+2}=\proxtP(\xi^{k+1})+\tilde{E}^{k+1},
         \end{equation*}
    	where 
        \begin{align*}
        	\tilde{E}^{k+1}=\left[
        	\begin{array}{c}
        		E^{k+1}\\
        		\sqrt{C}E^{k+1}
        	\end{array}\right];
        \end{align*}}
      \item {
       \begin{equation*}
      	{\rm dist}\left(\0,S_\Phi^{k+1}\left(\xi^{k+2}\right)\right)\le\left\Vert P^{-1}\right\Vert
      	{\rm dist}\left(\0,S_T^{k+1}\left(Z^{k+2}\right)\right);
       \end{equation*}
      }
      \item {
      	\begin{equation*}
      		\left\Vert\xi^{k+2}-\proxtP(\xi^{k+1})\right\Vert\le{\rm dist}\left(\0,S_\Phi^{k+1}\left(\xi^{k+2}\right)\right).
      	\end{equation*}
      	where
      	\begin{equation}\label{eq_def_sphi}
      		S_\Phi^{k+1}(\xi)=P^{-1}\Phi(\xi)+\left(\xi-P^{-1}Q\xi^{k+1}\right).
      	\end{equation}
      }
    \end{enumerate}
    \end{propo}  
    The proof is given in Appendix.  From the definition of $\tilde{E}^{k+1}$, we have
    \begin{equation*}
    	\left\Vert\tilde{E}^{k+1}\right\Vert=\left\Vert E^{k+1}\right\Vert_{I+C}\le
    	\rho\left(I+C\right)\left\Vert E^{k+1}\right\Vert.
    \end{equation*}
    Consequently, via Proposition \ref{propo_2}(a) and \ref{propo_2}(b), we obtain the following criteria of $\xi^{k+1}$ approximating $\proxtP(\xi^k)$, which corresponds to criteria (A1), (B1), (A2), and (B2), respectively,
    \begin{align*}
      &({\rm A3})\hspace{3em}
      \left\Vert\xi^{k+2}-\proxtP(\xi^{k+1})\right\Vert\le\rho(I+C)\varepsilon_{k+1},\quad
      \sum_{k=0}^\infty\varepsilon_{k+1}<\infty;\\
      &({\rm B3})\hspace{3em}
      \left\Vert\xi^{k+2}-\proxtP(\xi^{k+1})\right\Vert\le\rho(I+C)\delta_{k+1}
      \left\Vert\xi^{k+2}-\xi^{k+1}\right\Vert,
      \quad\sum_{k=0}^\infty\delta_{k+1}<\infty;\\
      &({\rm A4})\hspace{3em}
      {\rm dist}\left(\0,S_\Phi^{k+1}\left(\xi^{k+2}\right)\right)\le\left\Vert P^{-1}\right\Vert
      \varepsilon_{k+1}/\alpha,\quad\sum_{k=0}^\infty\varepsilon_{k+1}<\infty;\\
      &({\rm B4})\hspace{3em}
      {\rm dist}\left(\0,S_\Phi^{k+1}\left(\xi^{k+2}\right)\right)\le\left\Vert P^{-1}\right\Vert
      \left(\delta_{k+1}/\alpha\right)\left\Vert\xi^{k+2}-\xi^{k+1}\right\Vert,\quad\sum_{k=0}^\infty\delta_{k+1}<
      \infty.
    \end{align*}
    To be clear, the following map helps to clarify the relations among these criteria,
    \begin{equation}\label{map_A}
    \begin{aligned}
    	\begin{array}{rcr}
    		{\rm (A2)}& \stackrel{{\rm Eq.}\eqref{ine_estm}}{\Longrightarrow}& {\rm (A1)}\\
    		 {\rm Prop.}\ref{propo_2}(b)\Downarrow & ~ &{\rm Prop.}\ref{propo_2}(a)\Downarrow\\
    	    {\rm (A4)}& \stackrel{{\rm Prop.}\ref{propo_2}(c)}{\Longrightarrow}& {\rm (A3)}
    	\end{array},
    \end{aligned}
    \end{equation}
    \begin{equation}\label{map_B}
    \begin{aligned}
        \begin{array}{rcr}
        	{\rm (B2)}& \stackrel{{\rm Eq.}\eqref{ine_estm}}{\Longrightarrow}& {\rm (B1)}\\
        	{\rm Prop.}\ref{propo_2}(b)\Downarrow & ~ &{\rm Prop.}\ref{propo_2}(a)\Downarrow\\
        	{\rm (B4)}& \stackrel{{\rm Prop.}\ref{propo_2}(c)}{\Longrightarrow}& {\rm (B3)}
        \end{array}.
    \end{aligned}
    \end{equation}
    
     We shall also make use of the mapping
    \begin{equation*}
    	\Psi\triangleq I-{\rm prox}_{P^{-1}\Phi},
    \end{equation*}
    the following proposition provides some properties of $\Psi$ that are necessary for the convergence analysis.
    \begin{propo}[Proposition 1 \cite{Rockafellar76}]\label{propo_3}
    	\hspace{1em}
    	\begin{enumerate}[(a)]
    		\item {$\xi={\rm prox}_{P^{-1}\Phi}(\xi)+\Psi(\xi)$, and $\Psi(\xi)\in P^{-1}\Phi\left({\rm prox}_{P^{-1}\Phi}(\xi)\right)$ for all $\xi\in\mR^{N\times n}\times \mR^{N\times n}$.}\label{propo_3a}
    		\item {$\left\langle{\rm prox}_{P^{-1}\Phi}(\xi)-{\rm prox}_{P^{-1}\Phi}(\xi^\prime),
    			\Psi(\xi)-\Psi(\xi^\prime)\right\rangle\ge0$ for all $\xi, \xi^\prime$.}\label{propo_3b}
    		\item {$\left\Vert{\rm prox}_{P^{-1}\Phi}(\xi)-{\rm prox}_{P^{-1}\Phi}(\xi^\prime)\right\Vert^2+\left\Vert\Psi(\xi)-\Psi(\xi^\prime)\right\Vert^2\le\left\Vert\xi-\xi^\prime\right\Vert^2$ for all $\xi, \xi^\prime$.}\label{propo_3c}
    	\end{enumerate}
    \end{propo}

     The following theorem establishes the convergence to a solution to problem \eqref{P}, when the proximal steps of Proximal-Correction Algorithm are inexactly calculated under criterion (A1) (or (A2)). 
    \begin{thm}\label{thm_3}
    	Under Assumptions \ref{ass_1} and \ref{ass_2}, for any $\alpha>0$, let $\{z_i^k\}_{k\ge0}$ and $\{v_i^k\}_{k\ge0}$, $i=1, 2,\dots, N$, be any sequences generated by the inexact Proximal-Correction Algorithm under criterion $(\rm A1)$ \emph{(or $({\rm A2})$)}, then $\{z_i^k\}_{k\ge0}, i=1, 2,\dots, N$, converge to the same solution $z^\infty$ to problem \eqref{P}, and for each $i=1,2,\dots,N$, the sequence $\{v_i^k\}_{k\ge0}$ converges to one element of $T_i(z^\infty)$.
    \end{thm}
    {\rm\bf Proof} \ Checking the relations among criteria (A1)--(A4) which are shown in \eqref{map_A}, in fact, we only need to analyze the convergence under criterion (A3), and the convergence
    under other criteria is followed. The idea of this proof is along the lines of Theorem \ref{thm_2}.
    
    {\bf Step 1:} Let $z^*$ be an arbitrary solution to problem \eqref{P}, write $Z^*=\1\otimes(z^*)^\top$. From Lemma \ref{lem_1}, there exists $Y^*\in\mR^{N\times n}$ such that $\0\in\Phi(Z^*,Y^*)$. Denote $\xi^*=\left((Z^*)^\top,(Y^*)^\top\right)^\top$. Recall that $P\xi^*=Q\xi^*=\xi^*$, hence it holds that
    $Q\xi^*\in P\xi^*+\Phi(\xi^*)$, which implies 
    \begin{align*}
    \xi^*=\proxtP(\xi^*) \quad\text{and}\quad \Psi(P^{-1}Q\xi^*)=\Psi(\xi^*)=\0. 
    \end{align*}
    Plugging $\proxtP(\xi^{k+1})$ into
    $\Vert\xi^{k+1}-\xi^*\Vert$ and by criterion (A3), we have
    \begin{align}
    	\left\Vert\xi^{k+1}-\xi^*\right\Vert&\le\left\Vert\xi^{k+1}-\proxtP(\xi^k)\right\Vert+
    	\left\Vert\proxtP(\xi^k)-\proxtP(\xi^*)\right\Vert\notag\\
    	&\le\rho(I+C)\varepsilon_k+\left\Vert\xi^k-\xi^*\right\Vert,\quad (P\succcurlyeq Q\succ\0)\label{ineq_conv}
    \end{align}
    Summing \eqref{ineq_conv} over $k$ from $r$ to $l$ gives
    \begin{equation*}
    	\left\Vert\xi^{l+1}-\xi^*\right\Vert\le\left\Vert\xi^r-\xi^*\right\Vert+
    	\rho\left(I+C\right)\sum_{k=r}^l\varepsilon_k \quad\text{for all}\quad l>r.
    \end{equation*}
    Since $\sum_{k=0}^\infty\varepsilon_k<\infty$, it follows that the limit
    \begin{equation*}
    	\lim_{k\to\infty}\left\Vert\xi^k-\xi^*\right\Vert
    \end{equation*}  does exist, and
    \begin{equation*}
    		\left\Vert\xi^{l+1}-\xi^*\right\Vert\le
    		\left(\left\Vert\xi^0-\xi^*\right\Vert+\rho(I+C)\sum_{k=0}^\infty
    		\varepsilon_k\right)\triangleq M, \quad \forall l\ge0,
    \end{equation*}
    thus the sequence $\{\xi^k\}_{k\ge0}$ is bounded.
    Applying Proposition \ref{propo_3}\eqref{propo_3c} to $\xi=P^{-1}Q\xi^k$ and $\xi^\prime=P^{-1}Q\xi^*$ gets
    \begin{align*}
    	&\left\Vert\proxtP(\xi^k)-\xi^*\right\Vert^2+\left\Vert\Psi(P^{-1}Q\xi^k)\right\Vert^2\\
    	\le&\left\Vert P^{-1}Q(\xi^k-\xi^*)\right\Vert^2\\
    	\le&\Vert\xi^k-\xi^*\Vert^2.
    \end{align*}
    Hence, 
    \begin{align*}
    	&\left\Vert\Psi\left(P^{-1}Q\xi^k\right)\right\Vert^2-\left\Vert\xi^k-\xi^*\right\Vert^2+\left\Vert\xi^{k+1}
    	-\xi^*\right\Vert^2\\
    	\le&\left\Vert\xi^{k+1}-\xi^*\right\Vert^2-\left\Vert\proxtP(\xi^k)-\xi^*\right\Vert^2\\
    	=&\left\langle\xi^{k+1}-\proxtP(\xi^k),\left(\xi^{k+1}-\xi^*\right)+\left(\proxtP(\xi^k)-\xi^*\right)\right\rangle\\
    	\le&\left\Vert\xi^{k+1}-\proxtP(\xi^k)\right\Vert\left(\Vert\xi^{k+1}-\xi^*\Vert+\Vert\xi^k-\xi^*\Vert\right)
    	\\
    	\le&2M\rho(I+C)\varepsilon_k.
    \end{align*}
    Consequently, it follows that
    \begin{equation}\label{ineq_psi}
    	\left\Vert\Psi\left(P^{-1}Q\xi^k\right)\right\Vert^2\le\Vert\xi^k-\xi^*\Vert^2-\Vert\xi^{k+1}-\xi^*\Vert^2+2M\rho(I+C)\varepsilon_k.
    \end{equation}
    Summing \eqref{ineq_psi} over $k$ from $0$ to $\infty$, we obtain 
    \begin{equation*}
    	\sum_{k=0}^\infty\left\Vert\Psi(P^{-1}Q\xi^k)\right\Vert^2\le\left(\Vert\xi^0-\xi^*\Vert^2
    	+2M\rho(I+C)\sum_{k=0}^\infty\varepsilon_k\right)<\infty,
    \end{equation*}
    which entails
    \begin{equation*}
    	\lim_{k\to\infty}\left\Vert\Psi(P^{-1}Q\xi^k)\right\Vert=0.
    \end{equation*}
    Applying Proposition \ref{propo_3}(a) to $\xi=P^{-1}Q\xi^k$ deduces that
    \begin{equation*}
     \Psi\left(P^{-1}Q\xi^k\right)\in P^{-1}\Phi\left(\proxtP(\xi^k)\right).
    \end{equation*}
    Since $P^{-1}$ is symmetric positive definite, it's easy to verify that $P^{-1}\Phi$ is also maximal monotone.
    For any $\xi$ and $w$ satisfying $w\in P^{-1}\Phi(\xi)$, the monotonicity of $ P^{-1}\Phi$ suggests 
    \begin{equation}\label{ineq_monoty}
    	\left\langle\xi-\proxtP(\xi^k),w-\Psi\left(P^{-1}Q\xi^k\right)\right\rangle\ge0.
    \end{equation}
    Since $\{\xi^k\}$ is bounded, suppose $\xi^\infty$ is a cluster point of $\{\xi^k\}$, and
    \begin{equation*}
    \lim_{k\to\infty}\left\Vert\xi^{k+1}-\proxtP(\xi^k)\right\Vert=0,
    \end{equation*}
    it is also a cluster of $\left\{\proxtP(\xi^k)\right\}_{k\ge0}$.
    Taking the limit on the left side of \eqref{ineq_monoty}, we have
    \begin{equation*}
    	\langle\xi-\xi^\infty, w-\0\rangle\ge0 \quad\text{for all}\quad \xi, w\quad\text{satisfying}\quad w\in P^{-1}\Phi(\xi),
    \end{equation*}
    which means that $\0\in P^{-1}\Phi\left(\xi^\infty\right)$, i.e. $\0\in\Phi(\xi^\infty)$.
    
    {\bf Step 2:} This step is to show that there cannot be more than one cluster point of $\{\xi^k\}$, which is the same argument as in Theorem \ref{thm_1} and not repeated here. 
    
    {\bf Step 3:} Therefore, we assert that $\{\xi^k\}$ converges to $\xi^\infty$ with $\0\in\Phi\left(\xi^\infty\right)$, then Lemma \ref{lem_1} tells us that 
    \begin{align*}
    	\xi^\infty=\left[
    	\begin{array}{c}
    		Z^\infty\\
    		Y^\infty
    	\end{array}\right]&,
    	\quad\text{and}\quad Z^\infty=\1\otimes (z^\infty)^\top\quad\text{with}\quad \0\in\sum_{i=1}^NT_i(z^\infty),\\
    	&\lim_{k\to\infty}\Vert z_i^k-z^\infty\Vert=0,\quad\text{for all}\quad i=1,2,\dots,N. 
    \end{align*}
    Reviewing Eq.\eqref{eq_VZY}, since $\{Z^k\}$ and $\{Y^k\}$ converge to $Z^\infty$ and $Y^\infty$, respectively, it is obvious that there exists $V^\infty\in\mR^{N\times n}$ such that
    \begin{equation*}
    	\lim_{k\to\infty}V^k=V^\infty, \ \text{and} \ \alpha V^\infty+\sqrt{C}Y^\infty=\tilde{W}Z^\infty-Z^\infty=\0.
    \end{equation*}
    With the same argument as in Theorem \ref{thm_1}, the assumption of $V^\infty\notin T(Z^\infty)$ would lead to a contradiction to $\0\in\Phi(Z^\infty,Y^\infty)$, thus it must hold that
    \begin{equation*}
    	\lim_{k\to\infty}v_i^k=v_i^\infty, \ \text{and} \ v_i^\infty\in T_i(z^\infty), \ \text{for all} \ i=1,2,\dots, N.
    \end{equation*}
    \qed
    \subsection*{Convergence Rate}
    \begin{thm}\label{thm_4}
    	Under Assumptions \ref{ass_1} and \ref{ass_2}, for any $\alpha>0$, let $\{z_i^k\}_{k\ge0}$ and $\{v_i^k\}_{k\ge0}$, $i=1, 2,\dots, N$, be any sequences generated by the inexact Proximal-Correction Algorithm under criterion $(\rm B1)$ \emph{(or $({\rm B2})$)}. Assume that $\Phi^{-1}$ is Lipschitz continuous at $\0$ with modulus $a\ge0$; let 
    	\begin{equation*}
    		\mu=\frac{\Vert P\Vert a}{\sqrt{\left(\Vert P\Vert a\right)^2+1}}<1.
    	\end{equation*}
        Then the sequence $\{\xi^k\}_{k\ge0}$ defined by \eqref{eq_def_pq}, converges to $\xi^\infty$, the unique solution to $\0\in\Phi(\xi)$. Moreover, there is an index $\bar{k}$ such that
        \begin{equation*}
        	\left\Vert\xi^{k+1}-\xi^\infty\right\Vert\le\theta_k\left\Vert\xi^k-\xi^\infty\right\Vert \quad 
        	\text{for all} \ k\ge\bar{k},
        \end{equation*} 
       where 
       \begin{align*}
       	1>\theta_k&\equiv\frac{\mu+\rho(I+C)\delta_k}{1-\rho(I+C)\delta_k} \quad \text{for all} \ k\ge\bar{k},\\ 
       	\theta_k&  \to\mu\in (0,1) \ (\text{as} \ k\to\infty).
       \end{align*}
       Consequently, for each $i=1,2,\dots,N$, the sequence $\{z_i^k\}_{k\ge0}$ converges to $z^\infty$, the unique solution to $\0\in\sum_{i=1}^N T_i(z)$, with a linear rate, and the sequence $\{v_i^k\}_{k\ge0}$ converges to one element of $T_i(z^\infty)$. 
    \end{thm}
    {\bf Proof} \ Checking the relations among criteria (B1)--(B4) which are shown in \eqref{map_B}, it can be found that we only need to establish the convergence rate of $\{\xi^k\}_{k\ge0}$ under criterion (B3).
    
    The sequence $\{\xi^k\}_{k\ge0}$ satisfies criterion (A3) for $\varepsilon_k=\delta_k\Vert \xi^{k+1}-\xi^k\Vert$, so the conclusions of Theorem \ref{thm_3} are valid, that is, the sequence $\{\xi^k\}_{k\ge0}$ converges to $\xi^\infty$, the unique solution to $\0\in\Phi(\xi)$ ($\Phi^{-1}$ is Lipschitz continuous at $\0$). Next we will be devoted to estimating the convergence rate.
    
    Applying Proposition \ref{propo_3}\eqref{propo_3c} to $\xi=P^{-1}Q\xi^k$ and $\xi^\prime=P^{-1}Q\xi^\infty=\xi^\infty$ gives
    \begin{align}\label{ine_xifty}
    	&\left\Vert\proxtP(\xi^k)-\xi^\infty\right\Vert^2+\left\Vert\Psi(P^{-1}Q\xi^k)\right\Vert^2\notag\\
    	\le&\left\Vert P^{-1}Q(\xi^k-\xi^\infty)\right\Vert^2\notag\\
    	\le&\Vert\xi^k-\xi^\infty\Vert^2.
    \end{align}
    Then we have 
    \begin{equation*}
    	\left\Vert\Psi\left(P^{-1}Q\xi^k\right)\right\Vert\le\Vert\xi^k-\xi^\infty\Vert\to0.
    \end{equation*}
    Choose $\tilde{k}$ such that
    \begin{equation*}
    	\left\Vert\Psi\left(P^{-1}Q\xi^k\right)\right\Vert\le\tau \quad\text{ for all} \ k\ge\tilde{k}.
    \end{equation*}
    Since $\Phi^{-1}$ is Lipschitz continuous at $\0$ with modulus $a$, it's obvious that $\left(P^{-1}\Phi\right)^{-1}$ is also Lipschitz continuous at $\0$ with modulus $\Vert P\Vert a$.
    By Proposition \ref{propo_3}\eqref{propo_3a}, we have
    \begin{align*}
     \proxtP(\xi^k)\in\left(P^{-1}\Phi\right)^{-1}\left(\Psi(P^{-1}Q\xi^k)\right).
     \end{align*}
    Therefore, the Lipschitz condition \eqref{ine_Lip} can be invoked for $\xi=\proxtP(\xi^k)$ and $w=\Psi(P^{-1}Q\xi^k)$ when $k$ is sufficiently large:
    \begin{align*}
    	\left\Vert\proxtP(\xi^k)-\xi^\infty\right\Vert\le
    	a\Vert P\Vert\left\Vert\Psi\left(P^{-1}Q\xi^k\right)\right\Vert
    	\quad\text{for all} \ k\ge\tilde{k},
    \end{align*}
    which via \eqref{ine_xifty} yields 
    \begin{equation*}
    	\left\Vert\proxtP(\xi^k)-\xi^\infty\right\Vert^2\le\mu^2\left\Vert\xi^k-\xi^\infty\right\Vert^2,
    	\quad\forall k\ge\tilde{k},
    \end{equation*}
    or in other words,
    \begin{equation}\label{ine_mu}
      \left\Vert\proxtP(\xi^k)-\xi^\infty\right\Vert\le\mu\left\Vert\xi^k-\xi^\infty\right\Vert,
      \quad\forall k\ge\tilde{k}.
    \end{equation}
    Under criterion (B3) we have
    \begin{align*}
         &\left\Vert\xi^{k+1}-\proxtP(\xi^k)\right\Vert\\
      \le&\rho(I+C)\delta_k\Vert\xi^{k+1}-\xi^k\Vert\\
      \le&\rho(I+C)\delta_k\left(\Vert\xi^{k+1}-\xi^\infty\Vert+\Vert\xi^k-\xi^\infty\Vert\right).
    \end{align*}
    Thus via \eqref{ine_mu},
    \begin{align*}
    	\left\Vert\xi^{k+1}-\xi^\infty\right\Vert&\le\left\Vert\xi^{k+1}-\proxtP(\xi^k)\right\Vert
    	+\left\Vert\proxtP(\xi^k)-\xi^\infty\right\Vert\\
    	&\le\rho(I+C)\delta_k\left(\Vert\xi^{k+1}-\xi^\infty\Vert+\Vert\xi^k-\xi^\infty\Vert\right)
    	+\mu\left\Vert\xi^k-\xi^\infty\right\Vert,
    \end{align*}
    which is exactly
    \begin{equation*}
    		\left\Vert\xi^{k+1}-\xi^\infty\right\Vert\le\frac{\mu+\rho(I+C)\delta_k}{1-\rho(I+C)\delta_k}
    		\left\Vert\xi^k-\xi^\infty\right\Vert,\quad\forall k\ge\tilde{k}.
    \end{equation*}
    Since $\delta_k\to0$ and $\mu<1$, there exists $\bar{k}\ge\tilde{k}$ such that
    \begin{equation*}
    	\frac{\mu+\rho(I+C)\delta_k}{1-\rho(I+C)\delta_k}<1\quad\text{for all} \ k\ge\bar{k}.
    \end{equation*}
    Up to here, all the claimed conclusions are proven.\qed
    \subsection*{Example Supporting the Lipschitz Continuity of $\Phi^{-1}$}
    In Theorem \ref{thm_4}, we assume that $\Phi^{-1}$ is Lipschitz continuous at $\0$. A natural question  is: what conditions on operators $T_i, i=1,2,\dots,N$, can lead to this assumption? This question requires further investigation, here we give a supporting example by making use of the following proposition.
    \begin{propo}[Proposition 1 \cite{Robinson1981}]\label{propo_4}
    	Let $S:\mR^n\rightrightarrows\mR^m$ be a polyhedral multifunction, if $S(x_0)$ is single-valued, then there exists a constant $a\ge0$ such that $S$ is Lipschitz continuous at $x_0$ with modulus $a$.
    \end{propo}
    
    Considering a distributed quadratic programming with coupled constraints:
    \begin{align*}
    	&\hspace{2em}\min_{x\in\mR^p} \sum_{i=1}^N\frac{1}{2}\left\Vert A_ix-b_i\right\Vert^2\\
    	& {\rm s.t.}\quad\sum_{i=1}^N\left(C_ix-d_i\right)\le0,\\
    	& \hspace{2em} l_i\le x\le u_i, \ i=1,2,\dots,N,
    \end{align*}
   where $C_i\in\mR^{q\times p}, d_i\in\mR^q$. Denote $\Omega_i=\left\{x| \ l_i\le x\le u_i\right\},i=1,2,\dots,N$. We can write its dual problem \eqref{DP} satisfying
   \begin{align*}
   	\partial L_i(x,\lambda)&=\left[
   	\begin{array}{c}
   		A_ix-b_i+C_i^\top\lambda\\
   		-C_ix+d_i\\
   	\end{array}
   	\right]
   	+\left[
   	\begin{array}{c}
   		N_{\Omega_i}(x)\\
   		N_{\mR^q_+}(\lambda)
   	\end{array}\right]\\
    &=\tilde{A}_iz-\tilde{b_i}+N_{\tilde{\Omega}_i}(z).
   \end{align*}
   where $\tilde{\Omega}_i=\Omega_i\times\mR^q_+$, and $z=[x^\top,\lambda^\top]^\top$. Let $T_i(z)=\partial L_i(z), i=1,2,\dots,N,$ and $\xi=[Z^\top,Y^\top]^\top$. In this case, it's easy to show that the operator $\Phi$ can be formulated as 
   \begin{align*}
   	\Phi(\xi)=A\xi-B+N_D(\xi),
   \end{align*}
   where $N_D(\xi)$ is the cartesian product of $N_{\tilde{\Omega}_i}(z)$, $i=1,2,\dots,N$, and the single point set $\{\0\}$. 
   
   Suppose that the considered strongly convex programming admits a unique solution, i.e. $\Phi^{-1}(0)$ is single-valued. Since the normal cones $N_{\Omega_i}$, $i=1,2,\dots,N$, and $N_{\mR^q_+}$ are polyhedrons, we have that $\Phi$ and its inverse $\Phi^{-1}$ are polyhedral multifunctions. Therefore, the Lipschitz continuity of $\Phi^{-1}$ at $\0$ is guaranteed by applying Proposition \ref{propo_4} to $S=\Phi^{-1}$.
    
    \section{Numerical Simulations}
\begin{figure}[thbp]
    	\centering
    	\includegraphics[scale=0.25]{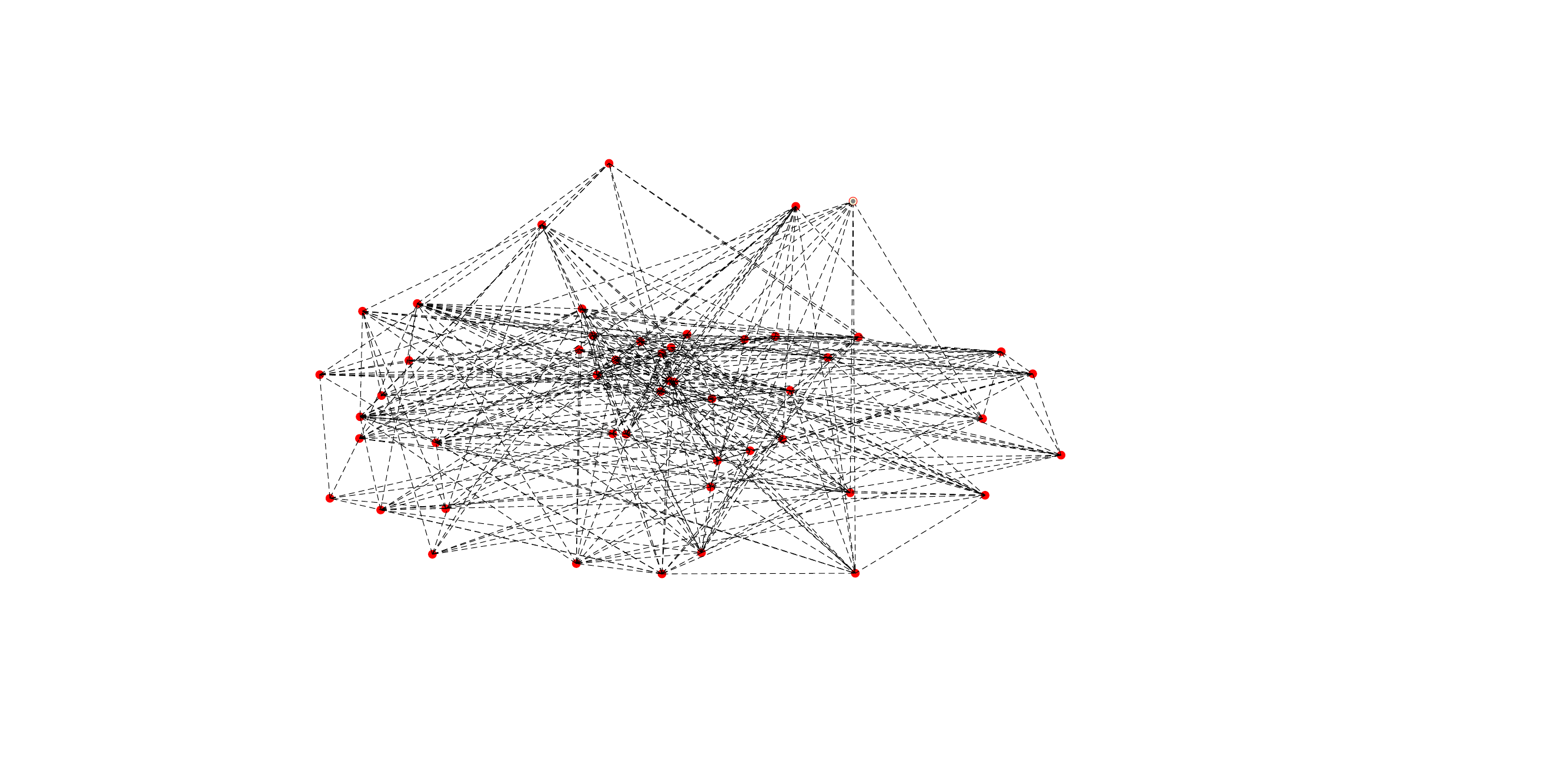}
    	\vspace{-3em}
    	\caption{Network topology of $50$ agents\qquad\quad}
    	\label{fig_1}
    \end{figure}

    In this section, we compare the practical performance of Proximal-Correction Algorithm with that of several alternative methods via the following distributed coupled constrained convex optimization problem: 
    \begin{equation}\label{eq_NS1}
    	\begin{aligned}
    	&\min_{x\in\mR}\sum_{i=1}^N a_i x\\
    	&{\rm s.t.}\quad \sum_{i=1}^N\left(-c_i\log(1+x)+\frac{b}{N}\right)\le0,\quad x\in\bigcap_{i=1}^N\left[\frac{i}{N},3-\frac{i}{N}\right].
    	\end{aligned}
    \end{equation}
    where $N=50, b=\frac{1}{2}N\log2, c_i=i/(N+1), a_i=i/N, i=1,\dots,N$. This example is motivated by applications in wireless networks \cite{Mateos17}, in which the coupled inequality constraints are to ensure the quality of service. The considered network topology is shown in Figure \ref{fig_1} and gives the mixing matrices $W$ and $\tilde{W}$ as follows
    \begin{align*}
    	w_{ij}&=\left\{
    	\begin{array}{cl}
    		\frac{1}{{\rm max}\{{\rm deg}(i),{\rm deg}(j)\}+1}\quad&\text{if}\quad j\neq i \quad\text{and}\quad (i,j)\in\E,\\
    		0,\quad&\text{if}\quad j\neq i\quad\text{and}\quad (i,j)\notin\E,\\
    		1-\sum_{l\in\V}w_{il} \quad&\text{if}\quad j=i.
    	\end{array}
    	\right.\\
    	\tilde{W}&=\frac{I+W}{2}.
    \end{align*}
    where ${\rm deg}(i)$ denotes the degree of agent $i$. 
    
    It is easy to verify that Slater's condition holds for \eqref{eq_NS1}. Then we need to solve a special case of \eqref{P}, i.e., $T_i=\partial L_i$,
    to obtain a solution to \eqref{eq_NS1}, where 
    \begin{equation*}
    	L_i(x,y)=a_ix+y\left(-c_i\log(1+x)+\frac{b}{N}\right)+I_{[i/N, 3-i/N]}(x)-I_{\mR_+}(y).
    \end{equation*}
    In this case, finding a pair $(u,v)\in\mR\times\mR$ such that $(u,v)=\left(I+\alpha\partial L_i\right)^{-1}(x,y)$ is equivalent to solving the following min-max subproblem
    \begin{equation*}
    	 (u,v)=\mathop{\arg\min\max}\limits_{\zeta\in\mR,  \eta\in\mR}
    	\left\{L_i(x,y)+\frac{1}{2\alpha}\Vert\zeta-x\Vert^2-\frac{1}{2\alpha}\Vert\eta-y\Vert^2\right\},
    \end{equation*}
     where the objective function is differentiable and strongly convex-concave, thus we can quickly obtain its unique solution by employing the primal-dual gradient method or other efficient algorithms.
     \subsection*{ Performance of Proximal-Correction Algorithm and Several Alternative Algorithms}
     \begin{figure}[htbp]
     	\centering
     	\subfigure[Evolution of agents' trial solutions $x_i^k$'s, ${\color{red}\triangleleft}$ denotes the optimal solution $x^*$.]{
     		\begin{minipage}[t]{0.35\linewidth}
     			\centering
     			\includegraphics[width=2\linewidth]{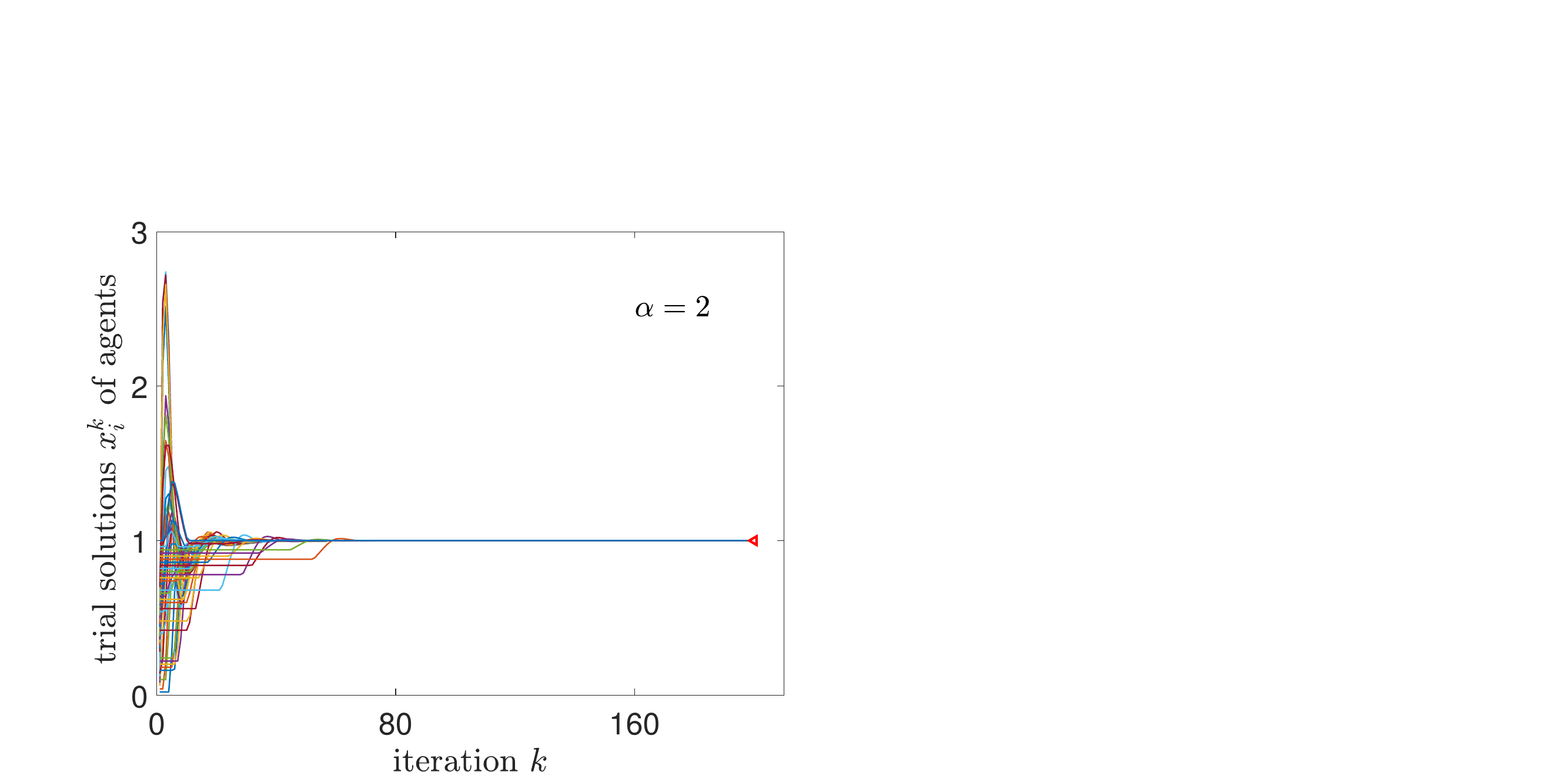}\label{fig_2a}
     			
     		\end{minipage}
     	}
     	\subfigure[Evolution of the constraint violation.]{
     		\begin{minipage}[t]{0.35\linewidth}
     			\centering
     			\includegraphics[width=2\linewidth]{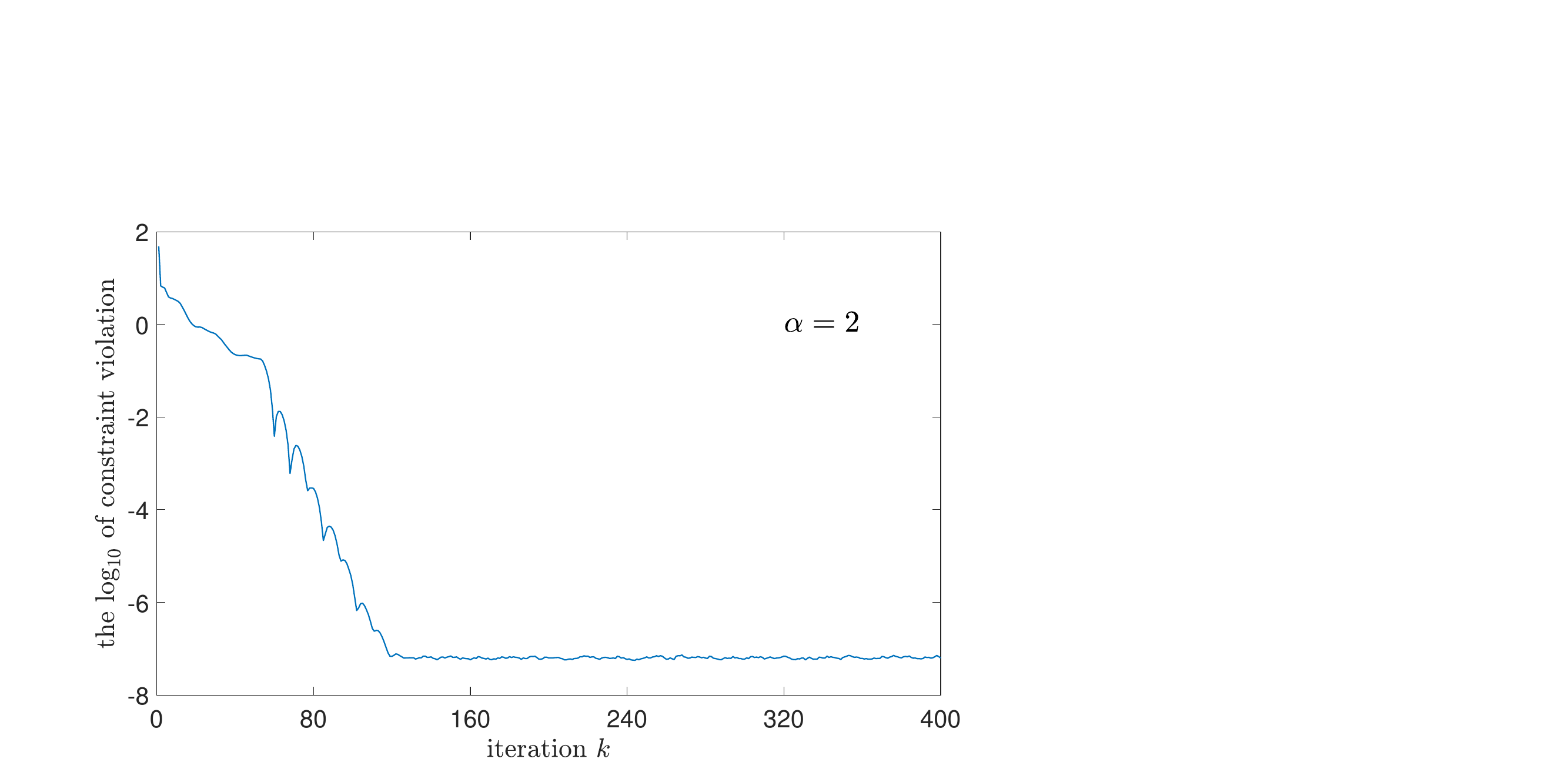}\label{fig_2b}
     		\end{minipage}
     	}
     	
     	\subfigure[Convergence performance under various step-sizes.]{
     		\begin{minipage}[t]{0.4\linewidth}
     			\centering
     			\includegraphics[width=1.7\linewidth]{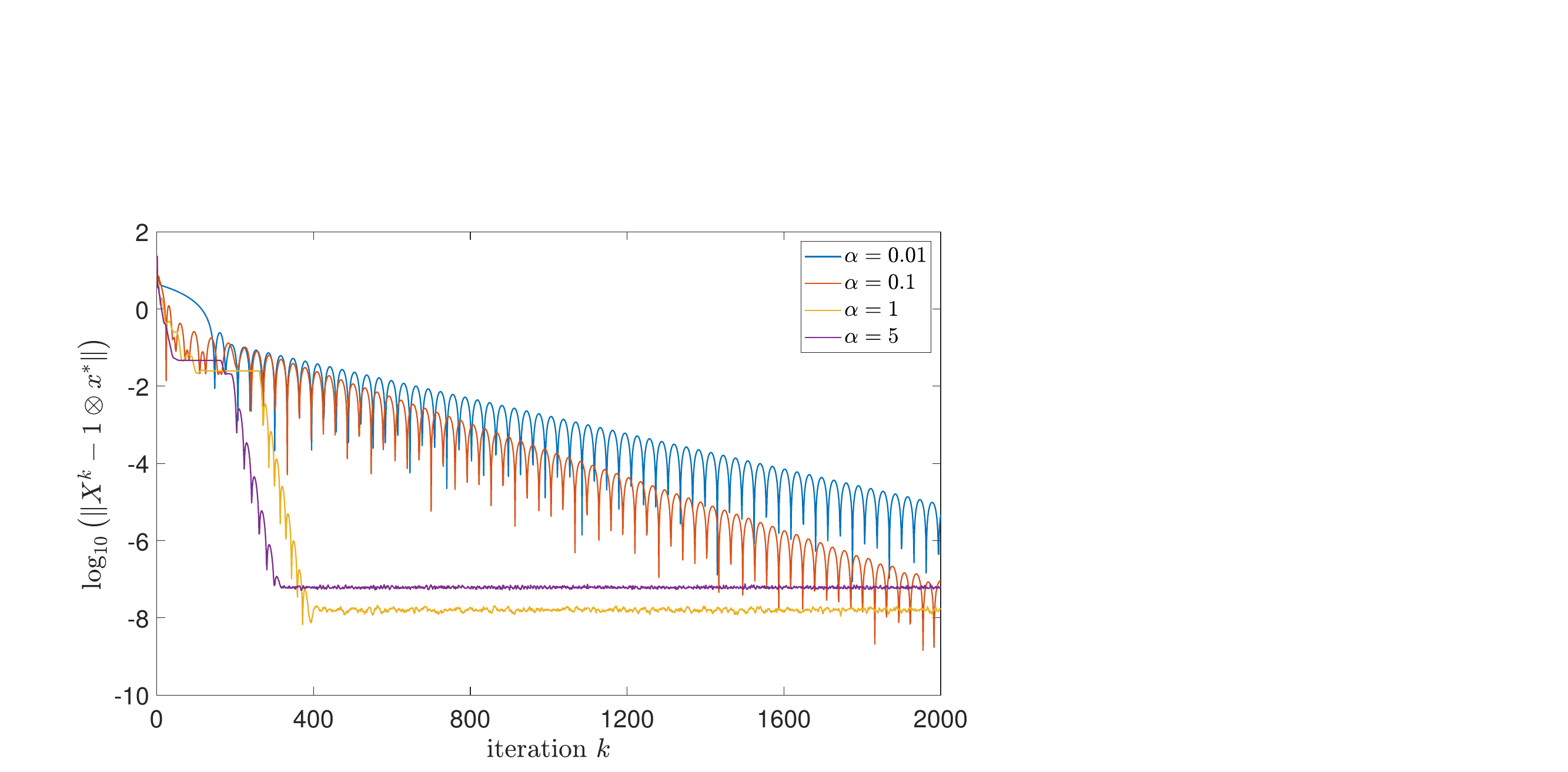}\label{fig_2c}
     		\end{minipage}
     	}
     	\subfigure[Evolution of constraint violation under various step-sizes.]{
     		\begin{minipage}[t]{0.4\linewidth}
     			\centering
     			\includegraphics[width=1.7\linewidth]{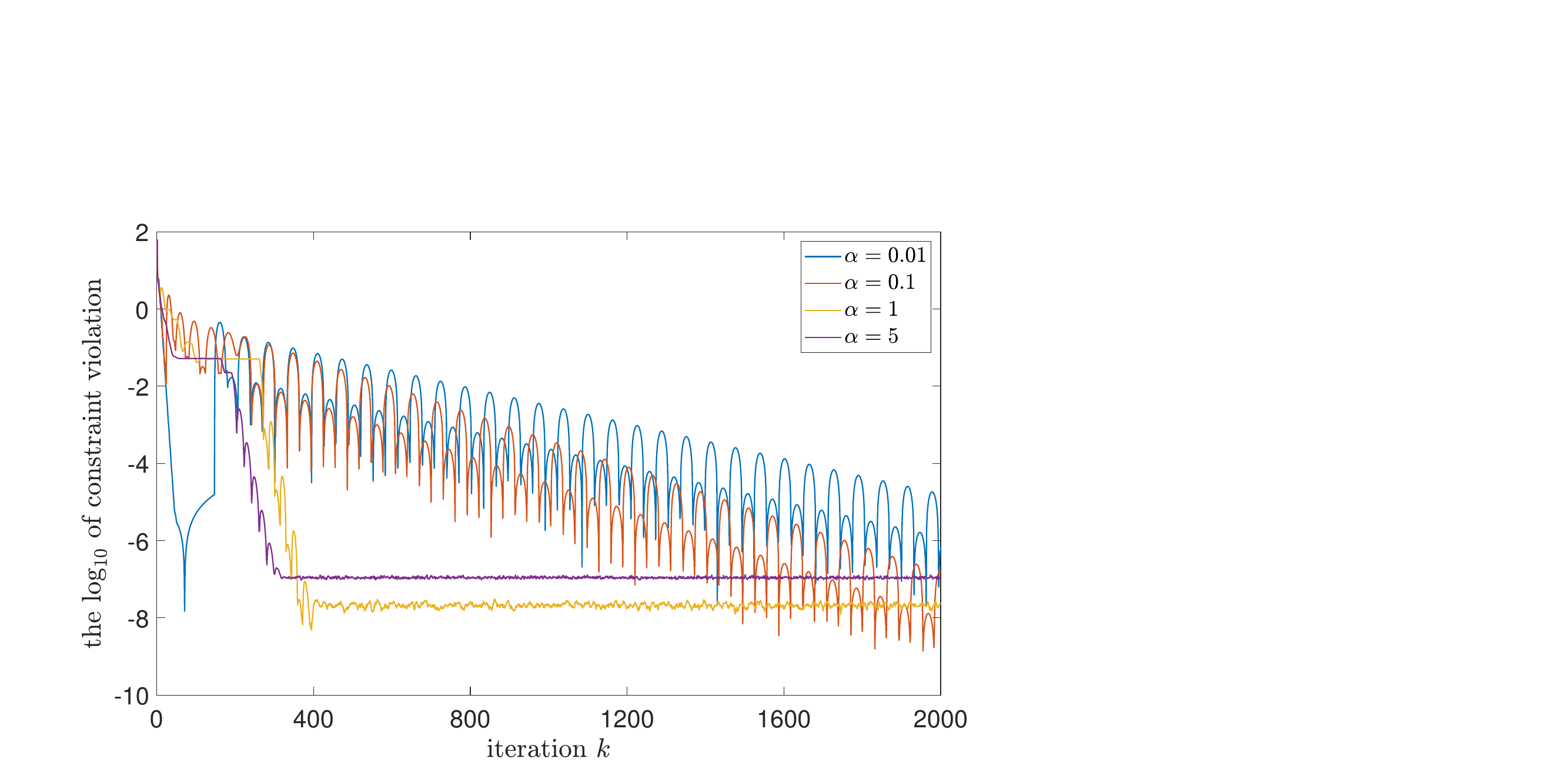}\label{fig_2d}
     		\end{minipage}
     	}
     	\centering
     	\caption{Convergence performance of Proximal-Correction}
     	\label{fig_2}
     \end{figure}
    In Fig. \ref{fig_2}, we first employ the Proximal-Correction algorithm with penalty parameter $\alpha=2$ to solve  \eqref{eq_NS1}. Figures \ref{fig_2a} and \ref{fig_2b} indicate that as information exchange proceeds among agents, our algorithm forces all agents to adjust their trial solutions towards the direction of being consensual, feasible, and optimal. The constraint violation is the sum of 
    $\Vert X^k-\1\otimes\left(\bar{x}^k\right)^\top\Vert$ and
    $\max\left(\sum_{i=1}^N(-c_i\log(1+x_i^k)+b/N),0\right)$,
    where $X^k=\left[x_1^k,\dots, x_N^k\right]^\top$ and $\bar{x}^k=\frac{1}{N}\sum_{i=1}^Nx_i^k$. Then we execute the Proximal-Correction algorithm with various values of the penalty parameter $\alpha$.
     
   Due to the absence of strong convexity of the objective functions, in Fig.\ref{fig_2c} and Fig.\ref{fig_2d}, the algorithm oscillates when the penalty parameter is  very small. As the penalty parameter gradually increases, the oscillations decrease and the proposed algorithm converges faster, which verifies that PRoximal-Correction converges for any $\alpha>0$ and inherits the features of the proximal point algorithm very well.
    
    References \cite{LiXX21}, \cite{ShuL20}, and \cite{WuXY22} also considered distributed convex optimization with coupling constraints. To employ the algorithms of \cite{ShuL20} and \cite{WuXY22} to problem \eqref{eq_NS1}, we can introduce a consensus equality constraint $(I-W)X=\0$, which means $x_1=x_2=\cdots=x_N$ and is also coupled, to reformulate \eqref{eq_NS1} as follows
    \begin{equation}\label{eq_NS2}
    	\begin{aligned}
    	&\min_{x_1, x_2,\dots, x_N}\sum_{i=1}^Na_i x_i\\
    	&{\rm s.t.}\quad\sum_{i=1}^N\left(-c_i\log(1+x_i)+\frac{b}{N}\right)\le0,\\
    	&\quad\quad\sum_{i=1}^Nl_i(x_i)=0,\quad x_i\in\left[\frac{i}{N},3-\frac{i}{N}\right],
    	\end{aligned}
    \end{equation}
    where 
    \begin{align*}
    	l_i(x_i)=\left[
    	\begin{array}{c}
    		L_{1i}x_i\\
    		\vdots\\
    		L_{Ni}x_i
    	\end{array}\right],\quad
    L=I-W.
    \end{align*}

    We execute Proximal-Correction, DPPD \cite{LiXX21}, the algorithm in \cite{ShuL20}, and IPLUX \cite{WuXY22} to solve \eqref{eq_NS1} or its equivalent reformulation \eqref{eq_NS2}. All the algorithm parameters are fine-tuned in their theoretical ranges to achieve the best possible convergence performance. Moreover, we let all the algorithms start from the same initial primal iterates for fairness.

\begin{figure}[htbp]
    	\centering
    	\subfigure[Evolution of solution errors of various algorithms]{
    			\begin{minipage}[t]{0.4\linewidth}
    			\centering
    			\includegraphics[width=1.7\linewidth]{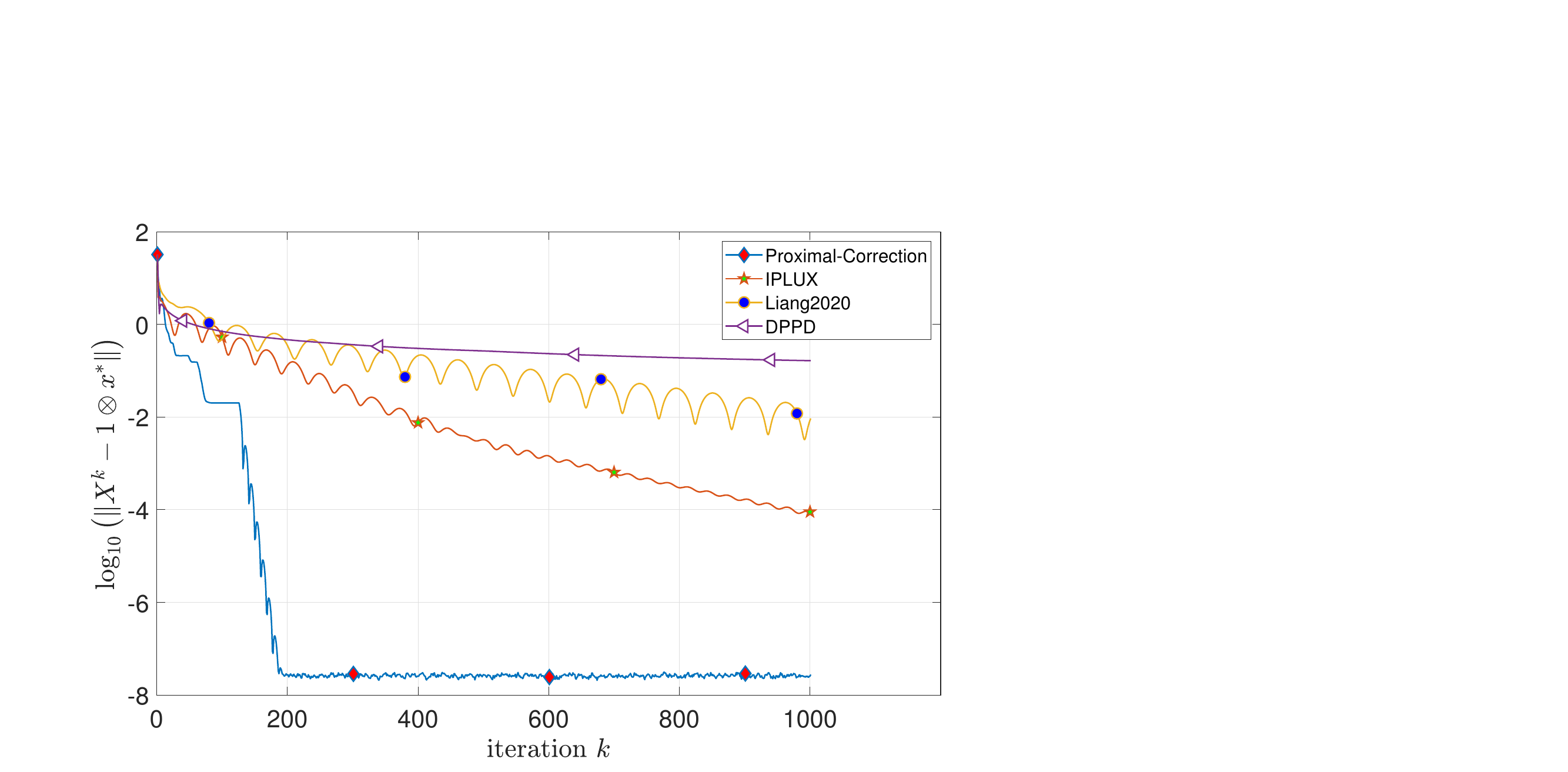}\label{fig_3a}
    		\end{minipage}
    	}
       \subfigure[Evolution of constraint violation of various algorithms]{
       	      \begin{minipage}[t]{0.4\linewidth}
       	      	\centering
       	      	\includegraphics[width=1.7\linewidth]{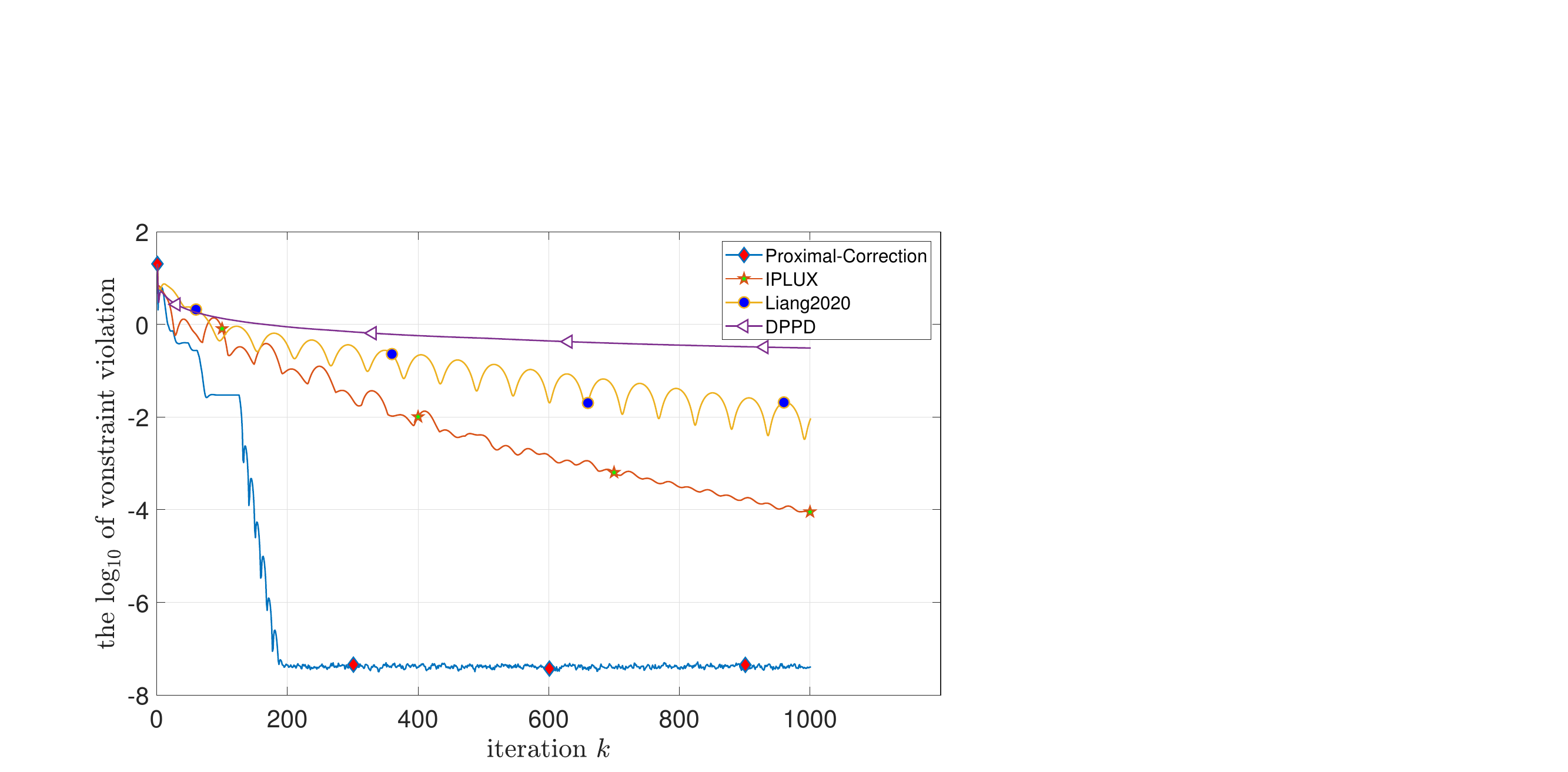}\label{fig_3b}
       	      \end{minipage}
       }
       \centering
       \caption{Comparison of the convergence performance of various algorithms}
       \label{fig_3}
    \end{figure}

    In Fig. \ref{fig_3}, we plot the solution errors and constraint violations generated by the aforementioned algorithms during 1000 iterations. As we can see, DPPD \cite{LiXX21} converges most slowly, because it adopts diminishing penalty parameters. The solution error and constraint violation of IPLUX \cite{WuXY22} both reach $10^{-4}$, which is matching its convergence rate of $O(1/k)$. While the proposed Proximal-Correction algorithm reaches an accuracy of $10^{-7}$  only needs about 200 iterations. 
    Compared to the three existing algorithms, Proximal-Correction presents a prominent advantage in convergence speed and accuracy concerning both optimality and feasibility.

    \subsection*{Performance of Inexact Proximal-Correction Algorithm}
    We also execute the inexact Proximal-Correction algorithm using criterion (A2) for problem \eqref{eq_NS1} to show its practical performance. A measure is provided in Remark \ref{rem_1} to verify condition (A2). For $i=1,2,\dots,N$, the precisions $\varepsilon_{k+1}^i$ of agent $i$ approximately calculating $\proxTi\left(\hat{z}_i^{k+1}\right)$ at iteration $k+1$, are simply set to be identical, i.e. $\varepsilon_{k+1}^i=\varepsilon_{k+1}^j\equiv\varepsilon_{k+1}$ for all $i,j\in\V$, where $\hat{z}_i^{k+1}$ is the $i$th row of $\widehat{Z}^{k+1}$. We adopt four precision sequences $\{\frac{1}{k+1}\}$, $\{\frac{1}{(k+1)^{1.5}}\}$, $\{\frac{1}{(k+1)^2}\}$, and $\{\frac{1}{(k+1)^{2.5}}\}$, for numerical simulations, and plot the solution errors and constraint violations to compare.
    
\begin{figure}[htbp]
    	\centering
    	\subfigure[Evolution of solution errors]{
    		\begin{minipage}[t]{0.4\linewidth}
    		\centering
    		\includegraphics[width=1.7\linewidth]{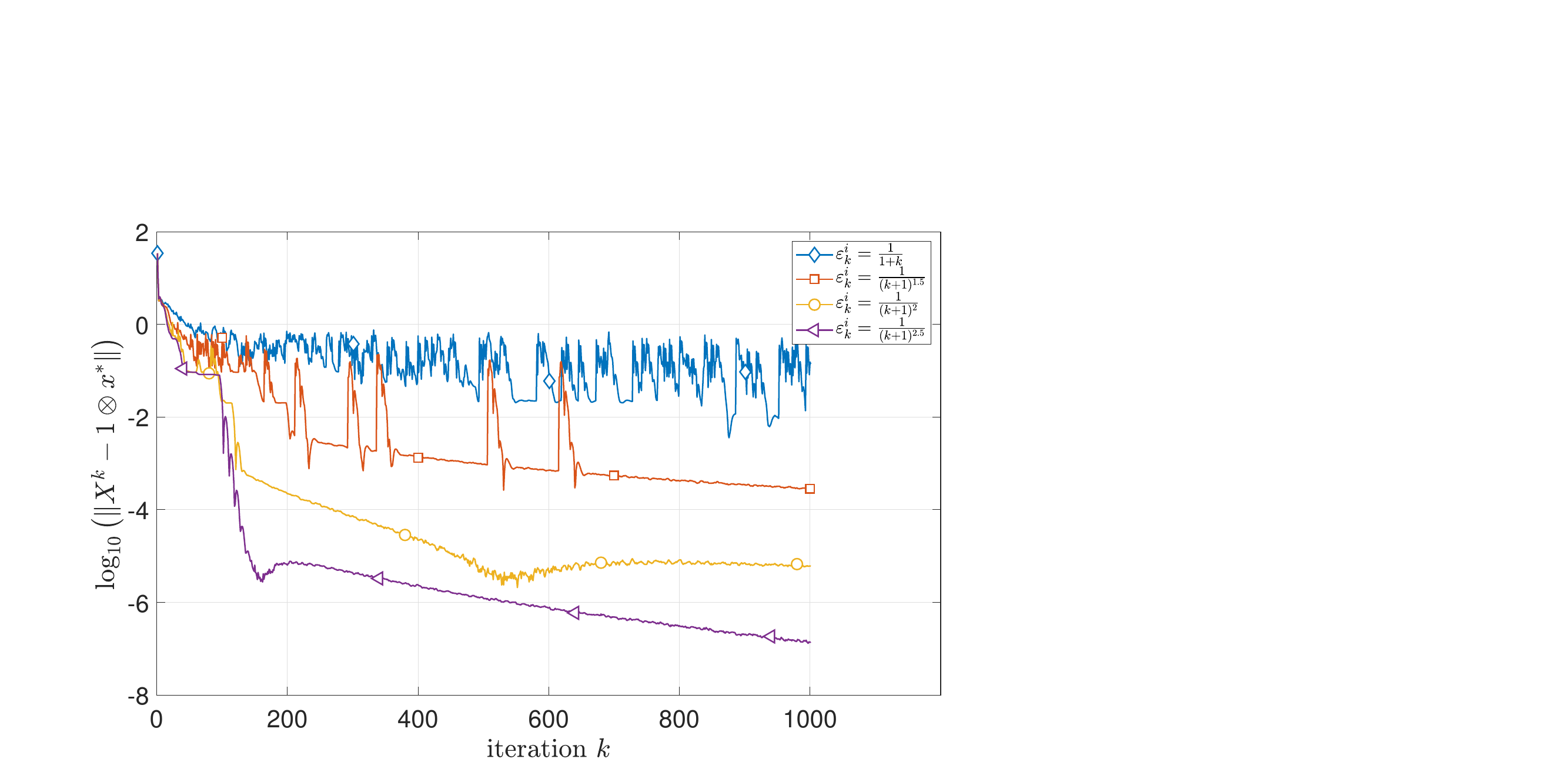}\label{fig_4a}
    	\end{minipage}
        }
       \subfigure[Evolution of constraint violation]{
       	\begin{minipage}[t]{0.4\linewidth}
       		\centering
       		\includegraphics[width=1.7\linewidth]{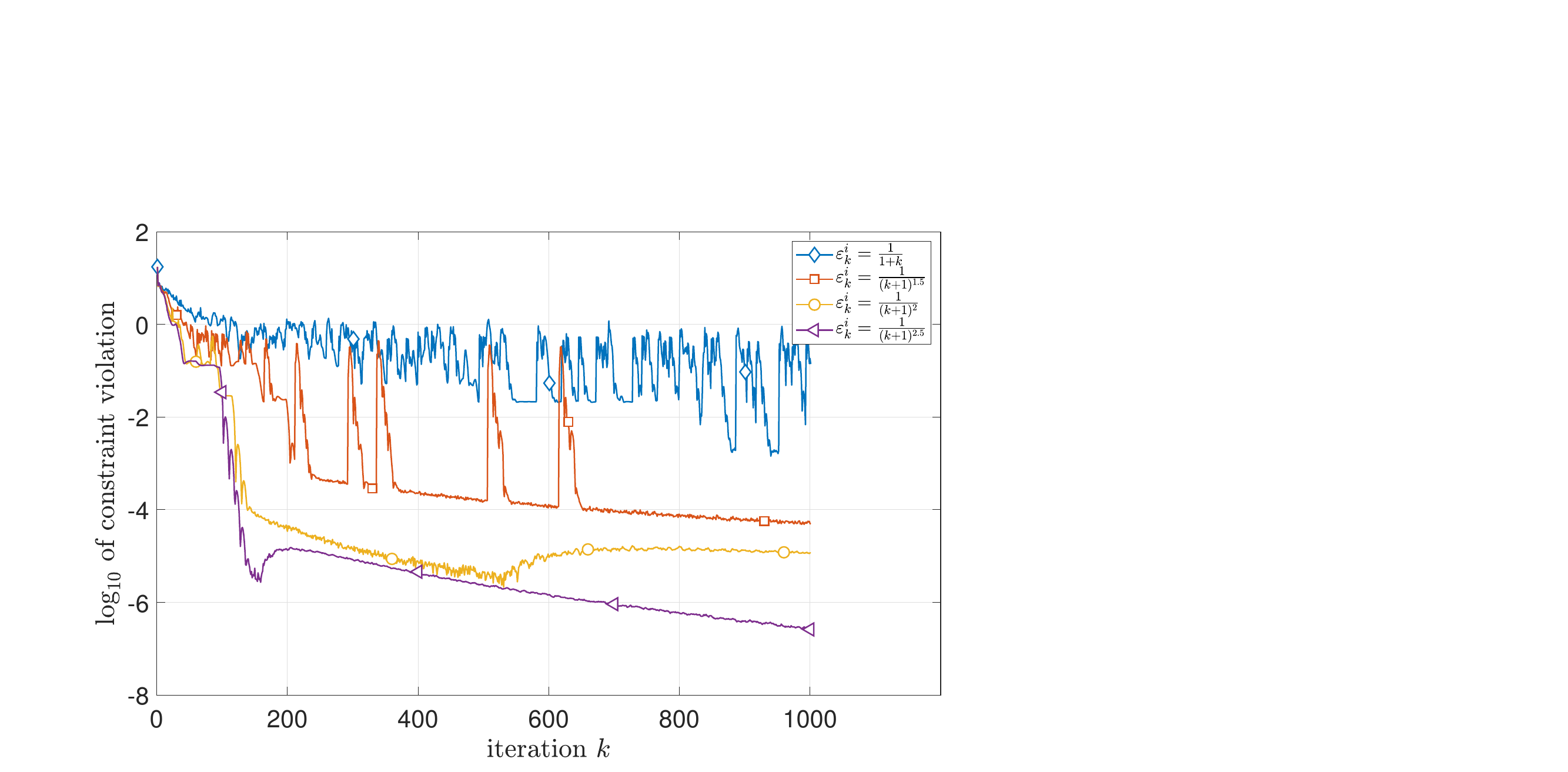}\label{fig_4b}
       	\end{minipage}
       }
       \centering
       \caption{Comparison of the convergence performance of the inexact Proximal-Correction under criterion (A2) with different precision sequences $\{\varepsilon_k\}_{k\ge0}$.}
       \label{fig_4}
    \end{figure}

    As shown in Fig.\ref{fig_4}, the convergence of the decision vector sequence $\{x_i^k\}_{k\ge0}$ failes when the precision sequence $\{\frac{1}{1+k}\}_{k\ge0}$ only has $\frac{1}{1+k}\to 0$ ( as $k\to\infty$), but $\sum_{k=0}^\infty \frac{1}{k+1}=\infty$. When the selected precision sequence $\{\varepsilon_k\}$ is summable, the convergence is guaranteed. The faster the precision sequence $\{\varepsilon_k\}_{k\ge0}$ converges to zero, i.e., the smaller the errors of computing the proximal mapping, the faster the optimality and feasibility of the decision vectors $\{x_i^k\}_{k\ge0}$, $i=1,2,\dots,N$, is to be satisfied. Fig.\ref{fig_4} also shows that the inexact Proximal-Correction algorithm is also very efficient, as we can see, when the precision sequence is $\{\frac{1}{(1+k)^2}\}$, the solution errors and constraint violation both reach an accuracy of $10^{-4}$ only needs about 400 iterations, which exhibits a faster speed compared with those three algorithms shown in Fig.\ref{fig_3}.

   \section{Conclusions}
    This paper considered finding a solution to the sum of maximal monotone operators in a multi-agent network, which can be linked to various types of distributed convex optimization problems. To this purpose, we developed the distributed Proximal-Correction algorithm, which extended PPA to the distributed settings and maintained fast convergence for any value of a constant penalty parameter. In addition, we established the convergence and the linear rate when the proximal steps of Proximal-Correction are inexactly calculated under two criteria. Although this paper considers an undirected graph, we would like to mention that Proximal-Correction can relax the doubly-stochasticity assumption on mixing matrices to singly-stochasticity by utilizing the push-sum protocol. Also, this relaxation requires further investigation. Our future research efforts would be to give general conditions that guarantee the Lipschitz continuity of $\Phi^{-1}$ at $\0$, and to modify Proximal-Correction to cope with a directed unbalanced network, even time-varying networks. Another research goal would be to relax the penalty parameter to be uncoordinated among agents. 
    \appendix
    \section{Proof of Proposition \ref{propo_2}}
    {\bf (a)} From \eqref{eq_zhat_y},
    \begin{align*}
    	Z^{k+2}&=\proxT\left(\widehat{Z}^{k+1}\right)+E^{k+1},\\
    	&=\proxT(\tilde{W}Z^{k+1}-\sqrt{C}Y^{k+1})+E^{k+1}.
    \end{align*}
    which follows 
    \begin{align*}
    	Q\xi^{k+1}=\left[
    	\begin{array}{c}
    		\tilde{W}Z^{k+1} \\
    		Y^{k+1}
    	\end{array}
    	\right]\in
    	\left[
    	\begin{array}{c}
    		Z^{k+2}-E^{k+1}+\alpha T(Z^{k+2}-E^{k+1})+\sqrt{C}Y^{k+1}\\
    		Y^{k+2}-\sqrt{C}Z^{k+2}
    	\end{array}
    	\right].
    \end{align*}
    Recall that $Y^{k+2}=Y^{k+1}+\sqrt{C}Z^{k+2}$, the term on the right side is equal to
    \begin{align*}
    	&\left[
    	\begin{array}{c}
    		(I-C)\left(Z^{k+2}-E^{k+1}\right)+\alpha T(Z^{k+2}-E^{k+1})+\sqrt{C}(Y^{k+2}-\sqrt{C}E^{k+1})\\
    		Y^{k+2}-\sqrt{C}E^{k+1}-\sqrt{C}\left(Z^{k+2}-E^{k+1}\right)
    	\end{array}
    	\right]\\
    	=&
    	\left[
    	\begin{array}{cc}
    		(I-C) & \0 \\
    		\0    & I
    	\end{array}
    	\right]
    	\left[
    	\begin{array}{c}
    		Z^{k+2}-E^{k+1}\\
    		Y^{k+2}-\sqrt{C}E^{k+1}
    	\end{array}\right]
    	+
    	\left[
    	\begin{array}{c}
    		\alpha T(Z^{k+2}-E^{k+1})+\sqrt{C}(Y^{k+2}-\sqrt{C}E^{k+1})\\
    		-\sqrt{C}\left(Z^{k+2}-E^{k+1}\right)
    	\end{array}
    	\right]\\
    	=&P\left(\xi^{k+2}-\tilde{E}^{k+1}\right)+\Phi\left(\xi^{k+2}-\tilde{E}^{k+1}\right),
    \end{align*}
    where 
    \begin{align*}
    	\tilde{E}^{k+1}=\left[
    	\begin{array}{c}
    		E^{k+1}\\
    		\sqrt{C}E^{k+1}
    	\end{array}\right].
    \end{align*}
    Hence it holds that $Q\xi^{k+1}\in(P+\Phi)\left(\xi^{k+2}-\tilde{E}^{k+1}\right)$, which implies 
    \begin{equation*}
    	\xi^{k+2}=\proxtP(\xi^{k+1})+\tilde{E}^{k+1}.
    \end{equation*}

    {\bf (b)} According to the definition of $S_\Phi^{k+1}$, we have
    \begin{align*}
    	S_\Phi^{k+1}\left(\xi^{k+2}\right)
    	&=P^{-1}\left(\Phi(\xi^{k+2})+P\xi^{k+2}-Q\xi^{k+1}\right)\\
    	&=P^{-1}\left[
    	\begin{array}{c}
    		\alpha T(Z^{k+2})+Z^{k+2}-\left(\tilde{W}Z^{k+1}-Y^{k+1}\right)\\
    		\0
    	\end{array}\right] \\
        &=P^{-1}\left[
        \begin{array}{c}
        	S_T^{k+1}\left(Z^{k+2}\right)\\
        	\0
        \end{array}\right]
    \end{align*}
    which gives
    \begin{equation*}
    	{\rm dist}\left(\0,S_\Phi^{k+1}\left(\xi^{k+2}\right)\right)\le\left\Vert P^{-1}\right\Vert
    	{\rm dist}\left(\0,S_T^{k+1}\left(Z^{k+2}\right)\right).
    \end{equation*}

    {\bf (c)} For any $\omega\in S_{P^{-1}\Phi}^{k+1}\left(\xi^{k+2}\right)$, it holds that
    \begin{equation*}
    	\omega+P^{-1}Q\xi^{k+1}\in\left(I+P^{-1}\Phi\right)\left(\xi^{k+2}\right),
    \end{equation*}
    and hence, 
    \begin{equation*}
    	\xi^{k+2}={\rm prox}_{P^{-1}\Phi}\left(w+P^{-1}Q\xi^{k+1}\right).
    \end{equation*}
    Then by virtue of the nonexpansiveness of ${\rm prox}_{P^{-1}\Phi}$ 
    \begin{align*}
    	&\left\Vert\xi^{k+2}-\proxtP(\xi^{k+1})\right\Vert\\
    	=&\left\Vert{\rm prox}_{P^{-1}\Phi}\left(\omega+P^{-1}Q\xi^{k+1}\right)-{\rm prox}_{P^{-1}\Phi}
    	\left(P^{-1}Q\xi^{k+1}\right)\right\Vert\\
    	\le&\Vert\omega\Vert
    \end{align*} 
    Thus
    \begin{equation*}
    	\left\Vert\xi^{k+2}-\proxtP(\xi^{k+1})\right\Vert
    	\le\min\left\{\Vert\omega\Vert: \omega\in S_{P^{-1}\Phi}^{k+1}\left(\xi^{k+2}\right)\right\}
    \end{equation*}
    as claimed.
    \bibliographystyle{unsrt}
    \bibliography{ref_DPCA}
\end{document}